\documentclass{amsart}
\usepackage{amssymb}
\usepackage{epsfig}

\title{Temperley-Lieb pfaffinants and Schur $Q$-positivity conjectures}
\date{December, 2006}
\author{Thomas Lam}
\email{tfylam (at) math (dot) harvard (dot) edu}
\thanks{T. L. was partially supported by NSF DMS--0600677.}

\author{Pavlo Pylyavskyy}
\email{pasha (at) math (dot) mit (dot) edu}
\date{}

\address{}


\thanks{}

\theoremstyle{plain}
\newtheorem{theorem}{Theorem}
\newtheorem{proposition}[theorem]{Proposition}
\newtheorem{prop}[theorem]{Proposition}
\newtheorem{lemma}[theorem]{Lemma}
\newtheorem{lem}[theorem]{Lemma}
\newtheorem{corollary}[theorem]{Corollary}
\newtheorem{conjecture}[theorem]{Conjecture}
\newtheorem{problem}[theorem]{Problem}

\theoremstyle{definition}
\newtheorem{definition}[theorem]{Definition}
\newtheorem{example}[theorem]{Example}

\theoremstyle{remark}
\newtheorem{remark}[theorem]{Remark}

\def\P{\mathbb P}
\def\N{\mathbb N}
\def\T{\mathbb T}
\def\wt{\mathrm{wt}}

\def\sh{\mathrm{sh}}

\def\X{\,\,\lower2pt\hbox{\input{figX.pstex_t}}}
\def\noXv{\,\,\lower2pt\hbox{
\begin{picture}(0,0)%
\includegraphics{figNoX.pstex}%
\end{picture}%
\setlength{\unitlength}{1973sp}%
\begingroup\makeatletter\ifx\SetFigFont\undefined%
\gdef\SetFigFont#1#2#3#4#5{%
  \reset@font\fontsize{#1}{#2pt}%
  \fontfamily{#3}\fontseries{#4}\fontshape{#5}%
  \selectfont}%
\fi\endgroup%
\begin{picture}(316,316)(293,-969)
\end{picture}
}}
\def\noXh{\,\,\lower2pt\hbox{
\begin{picture}(0,0)%
\includegraphics{figNoX2.pstex}%
\end{picture}%
\setlength{\unitlength}{1973sp}%
\begingroup\makeatletter\ifx\SetFigFont\undefined%
\gdef\SetFigFont#1#2#3#4#5{%
  \reset@font\fontsize{#1}{#2pt}%
  \fontfamily{#3}\fontseries{#4}\fontshape{#5}%
  \selectfont}%
\fi\endgroup%
\begin{picture}(316,316)(893,-369)
\end{picture}
}}
\def\noXDU{\,\,\lower2pt\hbox{
\begin{picture}(0,0)%
\includegraphics{figDU.pstex}%
\end{picture}%
\setlength{\unitlength}{1973sp}%
\begingroup\makeatletter\ifx\SetFigFont\undefined%
\gdef\SetFigFont#1#2#3#4#5{%
  \reset@font\fontsize{#1}{#2pt}%
  \fontfamily{#3}\fontseries{#4}\fontshape{#5}%
  \selectfont}%
\fi\endgroup%
\begin{picture}(316,316)(1493,-969)
\end{picture}
}}
\def\noXDD{\,\,\lower2pt\hbox{\input{figDD.pstex_t}}}
\def\noXUD{\,\,\lower2pt\hbox{
\begin{picture}(0,0)%
\includegraphics{figUD.pstex}%
\end{picture}%
\setlength{\unitlength}{1973sp}%
\begingroup\makeatletter\ifx\SetFigFont\undefined%
\gdef\SetFigFont#1#2#3#4#5{%
  \reset@font\fontsize{#1}{#2pt}%
  \fontfamily{#3}\fontseries{#4}\fontshape{#5}%
  \selectfont}%
\fi\endgroup%
\begin{picture}(316,316)(1493,-969)
\end{picture}
}}
\def\noXUU{\,\,\lower2pt\hbox{\input{figUU.pstex_t}}}
\def\noXRR{\,\,\lower2pt\hbox{\input{figRR.pstex_t}}}
\def\noXRL{\,\,\lower2pt\hbox{
\begin{picture}(0,0)%
\includegraphics{figRL.pstex}%
\end{picture}%
\setlength{\unitlength}{1973sp}%
\begingroup\makeatletter\ifx\SetFigFont\undefined%
\gdef\SetFigFont#1#2#3#4#5{%
  \reset@font\fontsize{#1}{#2pt}%
  \fontfamily{#3}\fontseries{#4}\fontshape{#5}%
  \selectfont}%
\fi\endgroup%
\begin{picture}(316,316)(1493,-969)
\end{picture}
}}
\def\noXLR{\,\,\lower2pt\hbox{
\begin{picture}(0,0)%
\includegraphics{figLR.pstex}%
\end{picture}%
\setlength{\unitlength}{1973sp}%
\begingroup\makeatletter\ifx\SetFigFont\undefined%
\gdef\SetFigFont#1#2#3#4#5{%
  \reset@font\fontsize{#1}{#2pt}%
  \fontfamily{#3}\fontseries{#4}\fontshape{#5}%
  \selectfont}%
\fi\endgroup%
\begin{picture}(316,316)(1493,-969)
\end{picture}
}}
\def\noXLL{\,\,\lower2pt\hbox{\input{figLL.pstex_t}}}
\def\XRR{\,\,\lower2pt\hbox{\input{figXRR.pstex_t}}}
\def\XRL{\,\,\lower2pt\hbox{\input{figXRL.pstex_t}}}
\def\XLR{\,\,\lower2pt\hbox{\input{figXLR.pstex_t}}}
\def\XLL{\,\,\lower2pt\hbox{\input{figXLL.pstex_t}}}

\def\P{{\bf \Pi}}
\def\Z{\mathbb{Z}}
\def\R{\mathbb{R}}

\def\type{{\rm type}}
\def\u{{\bf u}}
\def\p{{\bf p}}

\def\mult{{\rm mult}}
\def\cn{{\rm cn}}

\def\max{{\rm max}}
\def\uv{{\rm uv}}
\def\ph{{\rm ph}}
\def\sTL{\mathcal{T}}
\def\CP{P}
\def\Imm{\mathrm{Imm}}
\def\Pfaf{\mathrm{Pfaf}}
\def\TL{{\rm TL}}
\def\pf{\mathrm{pf}}
\def\D{\mathcal{D}}
\def\lex{{\rm lex}}
\def\tN{{\tilde{N}}}

\begin{document}

\begin{abstract}
We study pfaffian analogues of immanants, which we call pfaffinants.
Our main object is the TL-pfaffinants which are analogues of Rhoades
and Skandera's TL-immanants.  We show that $TL$-pfaffinants are
positive when applied to planar networks and explain how to
decompose products of complementary pfaffians in terms of
$TL$-pfaffinants.  We conjecture in addition that TL-pfaffinants have
positivity properties related to Schur $Q$-functions.
\end{abstract}

\maketitle
\section{Introduction}

An {\it {immanant}} of an $n \times n$ matrix $X = (x_{ij})$ is an
expression of the form
\begin{equation}\label{eq:imm}
\sum_{w \in S_n} f(w) x_{1,w(1)} \cdots x_{n,w(n)}
\end{equation}
where $f: S_n \longrightarrow \mathbb R$ is a function. The
well-known examples of immanants are determinants and permanents.
Desarmenien, Kung and Rota~\cite{DKR} gave a {\it {standard}} basis
of the space $I(X)$ of immanants, labeled by standard bitableaux
while recently Pylyavskyy~\cite{Pyl} introduced a basis labeled by
non-crossing bitableaux.


Immanants with certain positivity properties, most notably the
irreducible immanants, had been studied earlier
in~\cite{GJ,Gre,Hai,Ste92,SS}. In a series of papers \cite{Ska, RS1, RS2} Rhoades and Skandera studied the {\it {dual canonical basis}}
of $I(A)$, also called {\it Kazhdan-Lusztig immanants}, labeled by
permutations.  These immanants possess remarkable positivity
properties: (a) they are non-negative when applied to totally
non-negative matrices~\cite{RS1,RS2}, and (b) they are
Schur-positive when applied to Jacobi-Trudi matrices~\cite{RS2}.
This second property was used in~\cite{LPP} to resolve several
Schur-positivity conjectures.  The subset of the dual canonical
basis corresponding to $321$-avoiding permutations can be given a
purely combinatorial interpretation and were called {\it
{Temperley-Lieb immanants}}, or TL-immanants, in~\cite{RS1}. Rhoades
and Skandera also gave a simple positive combinatorial rule for
writing a product of two complementary minors of $A$ in terms
TL-immanants.

\medskip

The pfaffian $\pf(A)$ of a skew symmetric $2n \times 2n$ matrix $A$
(see Section~\ref{sec:pfaf}) replaces the symmetric group $S_{2n}$
in the determinant with the set of matchings of $2n$ points.
Replacing the symmetric group in (\ref{eq:imm}) with matchings one
also obtains a pfaffian analogue of immanants, which we call {\it
pfaffinants}. The main object of this paper are the TL-pfaffinants
denoted $\Pfaf_D(A)$, which are analogues of the TL-immanants.

Stembridge~\cite{St} interpreted the pfaffian $\pf(A(N))$ in terms
of non-intersecting path families in a planar network $N$, where
$A(N)$ is a skew-symmetric matrix obtained from $N$. Separately, it
is also known (\cite{JP,Mac}) that the Schur $Q$-function
$Q_\lambda$ is equal to the the pfaffian $\pf(A_\lambda)$ for a
particular skew symmetric matrix $A_\lambda$, which we call a {\it
$Q$-Jacobi-Trudi}- matrix.  Our search for the TL-pfaffinants
revolves around the following three properties:

\begin{enumerate}
\item
a product of complementary pfaffians should decompose positively and
simply in terms of the TL-pfaffinants;

\item
a TL-pfaffinant should be positive when evaluated on the skew
symmetric matrix $A(N)$ associated to a planar network;

\item
a TL-pfaffinant should be Schur $Q$-positive when evaluated on a
$Q$-Jacobi-Trudi matrix.
\end{enumerate}

The pfaffinants $\Pfaf_D(A)$ that we define satisfy properties (1)
and (2), and we conjecture that they satisfy property (3).  The
positivity properties (2) and (3) are subtly different from the
situation with TL-immanants.  Our definition of the pfaffinants
$\Pfaf_D(A)$ requires the intermediate definition of a {\it diagram
pfaffinant} $\Pfaf'_D(A)$.  It appears rather mysteriously that it
is the diagram pfaffinants that describe network and (conjecturally)
Schur $Q$-positivity.  We should point out that the correct pfaffian
analogue of the entire dual canonical basis is still missing.  A
basis of this entire space of pfaffinants (without the positivity
properties we desire) is given by DeConcini and Procesi~\cite{DP}
from the point of view of invariant theory.

One of the Schur $Q$-positivity conjectures
(Conjecture~\ref{conj:con2}) that we state is a Schur $Q$-function
version of a sequence of positivity results we call {\it cell
transfer}: the monomial positivity version is established
in~\cite{LP}, the fundamental quasi-symmetric function version
in~\cite{LP2} and the Schur positivity version in~\cite{LPP}.

\medskip

We now briefly describe the organization of the paper.  In Section
\ref{pfaf}, we define diagram pfaffinants and Temperley-Lieb
pfaffinants, and show that the latter form a basis for the space of
products of pairs of complementary pfaffians.  In Section
\ref{netw}, we explain Stembridge's work on pfaffians and planar
networks and show that TL-pfaffinants are non-negative when applied
to planar networks. We characterize the linear combinations of
products of pairs of complementary pfaffinants that are
network-nonnegative.  In Section \ref{rel} we explore the
relationship between TL-immanants and TL-pfaffinants when applied to
certain matrices. In Section \ref{pos} we state a number of
conjectures concerning Schur $Q$-positivity properties of
pfaffinants, and in addition we prove a number of intermediate
results.



\section{Pfaffians and Pfaffinants} \label{pfaf}

\subsection{Preliminaries}
\label{sec:pfaf}

A skew-symmetric matrix $A = (a_{ij})_{i,j=1}^n$ is a matrix
satisfying $A^t=-A$ or alternatively $a_{ij} = -a_{ji}$. These
matrices are in bijection with arrays $(a_{ij})_{1 \leq i < j \leq
n}$ obtained by taking the part of $A$ above the diagonal.  We
denote the corresponding array also by $A$ and will not usually
distinguish the skew-symmetric matrix from the upper-triangular
array.

Now suppose $A$ is a skew-symmetric $2n \times 2n$ matrix. Define
the {\it {pfaffian}} $\pf(A)$ of $A$ by
$$\pf(A) = \sum_{\pi \in F_{2n}} \epsilon(\pi) \prod_{(i,j) \in \pi}
a_{ij},$$ where the sum is taken over the set $F_{2n}$ of matchings
$\pi$ on $2n$ vertices, and $\epsilon(\pi)$ is the {\it sign} or
{\it crossing number} of a matching.  It can be determined by the
following rule: place the $2n$ vertices on a straight line and draw
all the edges in $\pi$ as arcs above this line.  Let ${\cn}(\pi)$
denote the number of crossings between the arcs. Then $\epsilon(\pi)
= (-1)^{{\cn}(\pi)}$. For convenience we write $a_\pi := \prod_{(i,j)
\in \pi}a_{ij}$ for any $\pi \in F_{2n}$.  We will generally think
of the matching $\pi$ as a set of unordered pairs of elements of
$[2n]$.  For example, if $n=2$ we have $\pf(A) =
a_{12}a_{34}-a_{13}a_{24}+a_{14}a_{23}$.

Let $I \subset [2n]$ be a $2m$-element subset and let $A_I$ be the
corresponding submatrix, obtained by taking only the rows and
columns with indices in $I$.  We denote by $\pf_I(A)$ the pfaffian
of this submatrix.  More generally, for disjoint subsets $I_1, I_2,
\ldots$ of $[2n]$ we denote $\pf_{I_1, I_2, \ldots}(A) =
\pf_{I_1}(A) \pf_{I_2}(A) \cdots$ the product of the corresponding
pfaffians.  If $I \subset [2n]$ we let $\bar I = [2n] \backslash I$
denote the complement of $I$ in $[2n]$.

Two special cases of $\pf_{I_1, I_2, \ldots}(A)$ are particular
important to us. One is the {\it complementary pfaffians} $\pf_{I,
\bar I}(A)$, which are the products of pfaffians of two {\it
{complementary}} subarrays. The second one is the {\it monomials}
$\pf_{\pi}(A) = a_\pi =  \prod_{(i,j) \in \pi} a_{ij}$.  Thus one
may also write the definition of the pfaffian as $\pf(A) = \sum_{\pi
\in F_{2n}} \epsilon(\pi) \pf_{\pi}(A)$.

Next, for an arbitrarily function $f: F_{2n} \to \R$ we define the
{\it {pfaffinant}} $\Pfaf_f(A) = \sum_{\pi \in F_{2n}} f(\pi)
\pf_{\pi}(A)$.  It is easy to see that if $I, J, \ldots $ is a
partitioning of $[2n]$ into disjoint sets (of even size), then
$\pf_{I,J,\ldots}$ is a pfaffinant. In particular, each $\pf_{\pi}$
is a pfaffinant with $f(\rho) = \delta_{\rho \pi}$.

\subsection{The complementary Pfaffian subspace}
Let $\R[A] = \R[a_{12},a_{13},\ldots]$ denote the $\R$-vector space
of polynomials in the variables $\{a_{ij}\}_{1 \leq i < j \leq 2n}$.
Now let $\CP_n:= \R[\pf_{I,\bar I }(A)] \subset \R[A]$ denote the
subspace spanned by the complementary pfaffians $\pf_{I, \bar I}$,
for all possible pairs $(I, \bar I)$, including the case $I =
\emptyset$.

We call a partitioning $(I, \bar I)$ of $[2n]$ {\it standard} if $I
= \{i_1 < i_2 < \cdots < i_a\}$ and $\bar I = \{j_1 < j_2 < \cdots <
j_b\}$ where $a \geq b$ and $i_k < j_k$ for each $k \in [1,a]$.
Alternatively, $(I, \bar I)$ is standard if $I$ and $\bar I$ form
the first and second rows of a standard Young tableau.  We say
$\pf_{I, \bar I}$ is standard if $(I, \bar I)$ is.

\begin{theorem}[\cite{DP}]
\label{thm:stan} A basis of $\CP_n$ is given by the set
$\{\pf_{I,\bar I}(A) \mid (I, \bar I) \; \text{is standard} \}$ of
standard complementary pfaffians. The dimension of $\CP_n$ over $\R$
is equal to the number of standard Young tableaux of size $2n$ with
at most $2$ rows, each row of even size.
\end{theorem}

\begin{proof}
In \cite{DP}, a product of several complementary pfaffians is
associated to any (possibly non-standard) tableau $T$ with even
parts. It is shown (\cite[Theorem 6.5]{DP}) that the set of such
products indexed by standard tableaux forms a basis for the space of
all pfaffinants. The straightening algorithm showing that any
tableau can be expressed in terms of standard ones (\cite[Lemma
6.1-6.3]{DP}) involves quadratic relations among products of
pfaffians.  Since the number of parts in the tableaux involved do
not increase in such straightenings, the statement of the theorem
follows.
\end{proof}

We will give another proof of Theorem \ref{thm:stan} later.

\begin{remark}
The following is the natural generalisation. Let $\CP_{k,n} \subset
\R[A]$ denote the subspace spanned by $k$ complementary pfaffians of
a $2n \times 2n$ skew-symmetric matrix. Then the dimension of
$\CP_{k,n}$ is equal to the number of standard tableaux of size $2n$
with at most $k$ rows such that each row has even length.
\end{remark}

\subsection{Symmetric Temperley-Lieb diagrams}

Consider a rectangle with the $2n$ points $1,2,\ldots,2n$ on the
left side and $2n$ points $1', 2', \ldots, 2n'$ on the right side
(the numbering goes from top to bottom). A {Temperley-Lieb diagram}
$D$ is a non-crossing matching on the resulting $4n$ vertices.  An
edge of $D$ is called {\it vertical} if it is of the form $(i,j)$ or
$(i',j')$ and is called {\it horizontal} if it is of the form $(i,
j')$.  A TL-diagram $D$ is {\it symmetric} if it has symmetry about
the vertical axis.  Thus all the horizontal edges in $D$ are of the
form $(i,i')$ and the vertical edges come in pairs
$\{(i,j),(i',j')\}$.  The {\it order} $|D|$ of a symmetric
TL-diagram is the number of edges in $D$ with both ends on the left
side of the rectangle or, alternatively, half the number of vertical
edges.  We call a TL-diagram $D$ {\it {even}} (or {\it odd})
depending on the order of $D$. We denote by $\sTL_n$ the set of
symmetric TL-diagrams on $4n$ vertices, and by $\sTL^e_n$ the subset
of even symmetric TL-diagrams.

\begin{proposition} \label{prop:size}
For any integer $n \geq 1$ we have $$|\sTL_n| = {{2n} \choose n}
\quad \text{and}  \quad |\sTL^e_n| = \frac{1}{2}{{2n} \choose n} =
{{2n-1} \choose n}.$$
\end{proposition}
\begin{proof}
We show that $\sTL_n$ is in bijection with $n$-subsets of a $2n$
element set. One possible such correspondence is obtained as
follows: for $D \in \sTL_n$ color all $i \in [2n]$ such that $(i<j)
\in D$ black.  Among the remaining points color black the largest
ones so that we get $n$ black points in total. The inverse map from
a coloring of $2n$ points black and white, $n$ of each color, can be
described as follows. Start reading the points in reverse order,
from $2n$ to $1$.  For each black point $i$ one encounters we find
the smallest $j > i$ colored white which has not yet been used and
include the edge $(i,j)$ in $D$.  If no such $j$ exists, we include
the edge $(i,i')$ in $D$.  After doing this for all the black
points, we include an edge $(j,j')$ for each unmatched white point
$j$.


Now let $\sTL^o_n = \sTL_n \backslash \sTL^e_n$ denote the set of
odd symmetric TL-diagrams. We define an involution $\omega$ on
$\sTL_n$ which sends $\sTL^e_n$ to $\sTL^o_n$. Let $D \in \sTL$.  If
$(1,1') \in D$, there exists some smallest $i \in [2n]$ where $i
\neq 1$ so that $(i, i') \in D$. We define $\omega(D)$ by removing
the edges $(1,1')$ and $(i, i')$ from $D$  and including the edges
$(1,i)$ and $(1', i')$. Otherwise $(1,k) \in D$ for some (even) $k
\in [2n]$. We define $\omega (D)$ by removing the edges $(1,k)$ and
$(1', k')$ and including the edges $(1,1')$ and $(k,k')$. The
involution $\omega$ shows that $|\sTL^e_n|=|\sTL^o_n| =
\frac{1}{2}|\sTL_n|$.
\end{proof}

\subsection{Diagram pfaffinants} \label{sect:dpfaf}

For each $D \in \sTL_n$ we now define a function $f_{D}: F_{2n}
\longrightarrow \mathbb Z$ which in turn gives us the {\it {diagram
pfaffinant}} $\Pfaf'_{D}(A) := \Pfaf_{f_D}(A)$.

Recall that we have $4n$ vertices on the sides of the rectangle: $1,
\ldots, 2n$ on the left side and $1', \ldots, 2n'$ on the right.
Given a matching $\pi \in F_{2n}$, let $\nu(\pi)$ be the matching on
$[2n] \cup [2n]'$ such that $(i, j')$ and $(i', j)$ are in
$\nu(\pi)$ if and only if $(i,j) \in \pi$. Pick a planar embedding
of $\nu(\pi)$ such that all edges lie inside the rectangle, and
every pair of edges intersect at most once.  We assume the embedding
is chosen (a) to have mirror symmetry, (b) no pair of edges have a
point of tangency, and (c) that no $3$ edges cross at a single
point. Call an embedding satisfying these conditions {\it {nice}}.
Such an embedding is far from unique, however we will show that the
construction does not depend on the choice of embedding. We assume
for now one such presentation has been chosen for each $\pi$, which
we will (abusing notation) denote by $\nu(\pi)$ as well.

The set of intersections among the edges of $\nu(\pi)$ can be
divided into two kinds: the {\it {unpaired crossings}}, which are
the crossings between pairs of edges of the form $(i, j')$ and $(i',
j)$; and the {\it {paired crossings}}, which are the pairs of
crossings between $(p',q)$ and $(r',s)$ and between $(p,q')$ and
$(r, s')$, where inequalities $q < s$ and $r < p$ either both fail or both hold.


Given $\pi \in F_{2n}$ we define a set $X(\pi)$ of uncrossings of
$\nu(\pi)$.  Each embedded graph $x \in X(\pi)$ is obtained from
$\nu(\pi)$ by uncrossing every intersection, where each intersection
can be uncrossed in two ways: as a vertical uncrossing ``\noXv'' or
as a horizontal uncrossing ``\noXh''.  In addition, we require that
paired crossings are uncrossed in the same way.  With this
additional restriction, the uncrossed diagram $x$ is mirror
symmetric.  Thus $x$ is topologically equivalent to an element $D(x)
\in \sTL_n$ union a number of closed loops.

We define the {\it {weight}} $\wt(x)$ of an uncrossed embedded graph
$x \in X(\pi)$ as
$$
\wt(x) = 2^{l(x)} (-1)^{\uv(x)+ \ph(x)}.
$$
Here $l(x)$ is the number of closed loops in $x$, where pairs of
mirror symmetric loops are counted only once; $\uv(x)$ is the number
of unpaired vertical uncrossings in $x$; and $\ph(x)$ is the number
of paired horizontal uncrossings in $x$.



Now we define $f_{D}: F_{2n} \to \Z$ by
$$
f_D(\pi) = \sum_{\substack{x \in X(\pi) \\ D(x) = D}} \wt(x).
$$

\begin{theorem}\label{thm:welldefined}
The function $f_{D}$ obtained in this way does not depend on the
particular embedding we have picked for each $\nu(\pi)$.
\end{theorem}

Theorem \ref{thm:welldefined} is in fact not logically required for
the rest of the paper.  Its proof is delayed to
Section~\ref{sec:Rei}.

\begin{figure}[ht]
\begin{center}
\input{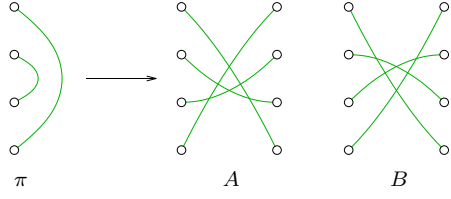}
\end{center}
\caption{Two different choices of the embedding $\nu(\pi)$ for $\pi
= \{(1,4),(2,3)\}$.} \label{fig:pfaf1}
\end{figure}

\begin{example}\label{ex:1423}
For $n = 2$ and $\pi = \{(1,4),(2,3)\}$, there are essentially two
different embeddings $A$ and $B$ of $\nu(\pi)$, shown in
Figure~\ref{fig:pfaf1}.  The embeddings are reflections of each
other about a horizontal axis.  These embeddings have two pairs of
mirror-symmetric crossings and two unpaired crossings, so the set
$X(\pi)$ has cardinality 16 in each case. The following table shows
the calculation of $f_D(\pi)$.

\begin{center}
\begin{tabular}{|c|c|c|}
\hline Diagram $D$ & $f_D(\pi)$ for embedding $A$ & $f_D(\pi)$ for
embedding $B$
\\
\hline ${\emptyset}$ & 1 & 1\\
\hline $\{(1, 2)\}$ & -1 & $2+2 - 1 -1 -2 - 1 = -1$ \\
\hline $\{(3, 4)\}$ & $2+2 - 1 -1 -2 - 1 = -1$ & -1\\
\hline $\{(1, 2), (3, 4)\}$ & $1 - 2 + 1 + 2 = 2$ &  $1 - 2 + 1 + 2 = 2$\\
\hline $\{(2, 3)\}$ &  $1+1-2 = 0$ & $1+1-2 = 0$\\
\hline $\{(2, 3), (1, 4)\}$ & -1  &-1\\
\hline
\end{tabular}
\end{center}
Thus for example for embedding $A$ there are 6 uncrossings $x \in
X(\pi)$ with $D(x) = \{(3,4)\}$.
\end{example}


\begin{example}\label{ex:diag}
For $n=2$ the diagram pfaffinants are given in the following table.
The diagrams are described by the sets of their vertical edges.  The
reader can verify that the coefficients of $a_{14}a_{23}$ agree with
the calculations in Example \ref{ex:1423}.

\begin{center}
\begin{tabular}{|c|c|}
\hline Diagram $D$ & Diagram pfaffinant $\Pfaf'_D(A)$
\\
\hline ${\emptyset}$ & $a_{12}a_{34}+a_{14}a_{23}-a_{13}a_{24}$ \\
\hline $\{(1, 2)\}$ & $-a_{14}a_{23}+a_{13}a_{24}-a_{12}a_{34}$ \\
\hline $\{(3, 4)\}$ & $-a_{14}a_{23}+a_{13}a_{24}-a_{12}a_{34}$ \\
\hline $\{(1, 2), (3, 4)\}$ & $a_{12}a_{34}+2a_{14}a_{23}-a_{13}a_{24}$ \\
\hline $\{(2, 3)\}$ & $0$ \\
\hline $\{(2, 3), (1, 4)\}$ & $a_{13}a_{24}-a_{14}a_{23}$ \\
\hline
\end{tabular}
\end{center}

\end{example}

We now state the main property of diagram pfaffinants.  Let $I
\subseteq [2n]$ and recall that $\bar I = [2n] \backslash I$ denotes
the complement of $I$ in $[2n]$.  The {\it $I$-coloring} of $[2n]
\cup [2n]'$ is obtained by coloring the elements of $I \cup {\bar
I}'$ black, and the elements $I' \cup \bar I$ white. We call a
diagram $D \in \sTL_n$ {\it compatible with $I$} (or simply
$I$-compatible) if each edge of $D$ has ends of different color in
the $I$-coloring.  Denote by $\D(I) \subset \sTL_n$ the set of
$I$-compatible diagrams.

\begin{theorem}
\label{thm:diag} Let $I \subset [2n]$ be a subset of even
cardinality. Then
$$
\pf_{I,\bar I}(A) = \sum_{D} \Pfaf'_{D}(A)$$
where the sum is over
all $I$-compatible diagrams of $\sTL_n$.
\end{theorem}
The following proof imitates a proof in~\cite{LPP}.
\begin{proof}
Let $\pi \in F_{2n}$.  Then the monomial $a_\pi$ occurs in
$\pf_{I,\bar I}$ if no edge of $a_\pi$ connects an element of $I$
with an element of $\bar I$.  In other words, $\pi$ must be the
union of the two matchings $\pi_I$ and $\pi_{\bar I}$ obtained by
restricting the vertex set.  The coefficient of $a_\pi$ in
$\pf_{I,\bar I}$ is then equal to $(-1)^{\cn(\pi_I) + \cn(\pi_{\bar
I})}$.

Now suppose $x \in X(\pi)$ is an uncrossing of $\nu(\pi)$ such that
$D(x) \in \D(I)$.  We direct all the strands and loops in $x$ so
that the initial vertex of each strand belongs to $I \cup (\bar I)'$
(and, thus the end vertex belongs to $\bar I \cup I'$).  We allow
the closed loops to be directed in either direction.  Thus the
coefficient of $a_\pi$ in $\sum_{D} \Pfaf'_{D}(A)$ is equal to the
sum of $(-1)^{\uv(y) + \ph(y)}$ over all orientations $y$ of the
uncrossings $\{x \in X(\pi) \mid D(x) \in \D(I)\}$.

\begin{figure}[ht]
\begin{center}
\input{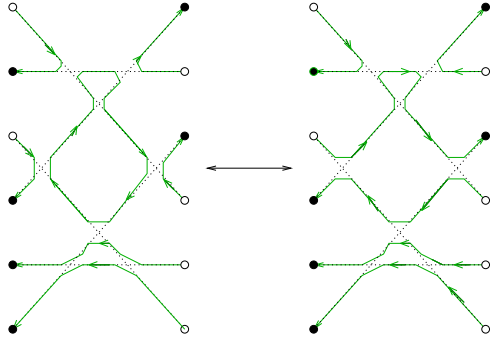}
\end{center}
\caption{The partial involution $\iota$ in the proof of
Theorem~\ref{thm:diag}.} \label{fig:pfaf2}
\end{figure}

Now we define a sign-reversing partial involution $\iota$ on this
set of oriented graphs. A {\it misaligned uncrossing\/} is an
uncrossing of the form ``\noXDU'', ``\noXUD'', ``\noXLR'', or
``\noXRL''. We say that we {\it {switch}} a misaligned uncrossing if
we apply one of the following transformations:
$\textrm{\noXDU}\longleftrightarrow \textrm{\noXRL}$ or
$\textrm{\noXUD}\longleftrightarrow \textrm{\noXLR}$. If $y$
contains any misaligned uncrossings then we let $\iota$ switch the
leftmost such uncrossing. If this uncrossing is a paired uncrossing,
we also switch its mirror image. If all the uncrossings are aligned,
then $\iota$ is not defined. Since $\iota$ is a sign-reversing
involution on the set of oriented graphs where it is defined, we
need only consider the contribution of $(-1)^{\uv(y) + \ph(y)}$ for
oriented graphs $y$ where $\iota$ is undefined. An example of the
application of $\iota$ for $n=3$ and $I = \{1,3\}$ is given in
Figure \ref{fig:pfaf2}. We switch the leftmost misaligned
uncrossing, which in this case happens to be paired.

Now suppose that $y_\pi$ is an oriented diagram with only aligned
uncrossings (see for example Figure \ref{pfaf3}). Then converting
the uncrossings back into crossings, keeping the orientation the
same, we obtain an orientation $\mu(\pi)$ of $\nu(\pi)$ such that
all edges start in $I$ end in $I'$ or start in $\bar I$ and end in
$(\bar I)'$.  Thus $\pi$ is the union of two matchings $\pi_I$ and
$\pi_{\bar I}$.  It is also clear that one can recover $y_\pi$ from
$\mu(\pi)$ and that $\mu(\pi)$ is completely determined by
$\nu(\pi)$.  Thus $y_\pi$, if it exists, is unique.

 \begin{figure}[ht]
\begin{center}
\input{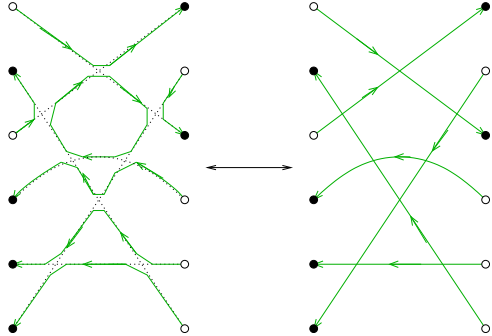}
\end{center}
\caption{Obtaining the orientation $\mu(\pi)$ of $\nu(\pi)$ from the
uncrossing $y_\pi$.} \label{pfaf3}
\end{figure}

Finally,  we calculate the sign of $y_\pi$.  Each unpaired crossing
of $\nu(\pi)$ corresponds to the intersection of $(i,j')$ with
$(i',j)$ for an edge $(i,j)$ in $\pi_I$ or $\pi_{\bar I}$.  These
crossings are always uncrossed horizontally to obtain $y_\pi$, and
so contributes no sign to $y_\pi$.  Each paired crossing $(c,c')$ in
$\nu(\pi)$ arises from a crossing $\xi$ of $\pi$.  To obtain
$y_\pi$, the pair $(c,c')$ is uncrossed horizontally if $\xi$ is a
crossing in $\pi_I$ or $\pi_{\bar I}$, and $(c,c')$ is uncrossed
vertically otherwise.  Thus $(-1)^{\uv(y_\pi) + \ph(y_\pi)} =
(-1)^{\cn(\pi_I) + \cn(\pi_{\bar I})}$, and we have checked that the
monomial $a_\pi$ appears in both sides with the same coefficient.
\end{proof}


\subsection{Temperley-Lieb pfaffinants}
Let $D \in \sTL_n$. For $i, j \in \{1, \ldots, 2n\}$ satisfying
$i<j$ we call the edge $(i,j)$ of $D$ {\it {odd}} if $i$ is odd and
{\it {even}} otherwise. For $D \in \sTL_n$ let $S(D)$ be the set of
all diagrams in $\sTL_n$ that can be obtained from $D$ by erasing
several odd edges (and their mirror images) and matching the
resulting unmatched vertices by horizontal edges of the form
$(i,i')$. In particular, $D \in S(D)$.

\begin{lem}
\label{lem:include} If $D_1, D_2 \in \sTL_n$ and $D_1 \in S(D_2)$
then $S(D_1) \subset S(D_2)$. The size of $S(D)$ is a power of $2$.
\end{lem}

\begin{proof}
The first statement is clear since after obtaining $D_1$ out of
$D_2$ by removing several odd edges, we can keep removing the
remaining odd edges, and the result belongs to $S(D_2)$ by
definition. For the second part, note that if $(i,j)$ is an odd
edge, that is if $i$ is odd, then all the edges inside $[i,j]$
cannot be removed either because they are even or because they are
contained within the segment bounded by ends of an even edge. Thus
all odd edges that can be removed can be removed independently one
from another, which implies the statement of the lemma.
\end{proof}

\begin{lem}
\label{lem:Imax} Suppose $D \in \sTL_n$ and $I \subset [2n]$ is a
subset of even cardinality.  If $D \in \D(I)$ then $D' \in \D(I)$
for every $D' \in S(D)$.  Conversely, if $D' \in \D(I)$ then there
exists a unique $D_\max \in \sTL^e_n \cap \D(I)$ such that $D' \in
S(D_\max)$ and $D_\max$ is maximal in the following sense: if $D'
\in S(D)$ for some other $D \in \D(I)$ then $S(D) \subset
S(D_\max)$.
\end{lem}
\begin{proof}
The first statement follows immediately from the definitions of the
set $S(D)$ and of $I$-compatibility.  Now let $D' \in \D(I)$. We say
that a vertex $i \in [2n]$ is free if $(i,i')$ is an edge in $D$. It
is clear that there are the same number of black and white vertices
in the $I$-coloring amongst the non-free vertices. Also, one checks
that the free vertices alternate in parity beginning with an odd
vertex and ending with an even vertex. If there are two vertices $i
< j$ such that between $i$ and $j$ there are no free vertices, $i$
is odd, $j$ is even and they have different colors then we call the
pair $(i,j) \in [2n] \times [2n]$ {\it addable}. Removing $(i,i')$
and $(j,j')$ from $D'$ and adding $(i,j)$ and $(i',j')$ gives some
$D \in \sTL_n \cap \D(I)$ such that $D' \in S(D)$.  The unique
maximal such $D = D_\max$ is obtained by performing the above
operation for every pair of addable vertices. Since $I$ is required
to have even cardinality and all the free vertices of $D_\max$ has
the same color, $D_\max$ must be even.
\end{proof}

We say that $D \in \sTL^e_n$ is {\it $I$-maximal} if it has the form
$D_\max$ as in Lemma~\ref{lem:Imax}.  We denote the set of
$I$-maximal diagrams by $\D_\max(I)$.  By Lemma~\ref{lem:include},
if $D_1, D_2 \in \D_\max(I)$ then $D_1 \notin S(D_2)$ and $D_2
\notin S(D_1)$.

\begin{definition}
Let $D \in \sTL^e_n$.  Define the {\it TL-pfaffinant} $\Pfaf_{D}(A)$
by
$$\Pfaf_{D}(A) = \sum_{D' \in S(D)} \Pfaf'_{D'}(A).$$
\end{definition}

\begin{example}\label{ex:even}
For $n=2$ the TL-pfaffinants are given in the following table,
calculated using Example~\ref{ex:diag}. The even diagrams are
described by the sets of their vertical edges.

\begin{center}
\begin{tabular}{|c|c|}
\hline Even diagram $D$ & TL-pfaffinant $\Pfaf_D(A)$
\\
\hline ${\emptyset}$ & $a_{12}a_{34}+a_{14}a_{23}-a_{13}a_{24}$ \\
\hline $\{(1, 2), (3, 4)\}$ & $a_{13}a_{24}-a_{12}a_{34}$ \\
\hline $\{(2, 3), (1, 4)\}$ & $a_{13}a_{24}-a_{14}a_{23}$ \\
\hline
\end{tabular}
\end{center}

\begin{example}
\label{ex:pfafhor} Let $I = [2n]$ and let $D \in \sTL^e_n$ be the
even symmetric TL-diagram with all edges horizontal.  Then
$\pf_{I,\bar I}(A) = \Pfaf_D(A)$.
\end{example}

\end{example}

\begin{theorem}
\label{thm:decomp} Suppose $I \subset [2n]$ is a subset with even
cardinality. Then
$$\pf_{I, \bar I}(A) = \sum_{D \in \D_\max(I)} \Pfaf_D(A).$$
\end{theorem}

\begin{proof}
By Theorem~\ref{thm:diag}, it suffices to show that the set of
$I$-compatible diagrams $\D(I) \subset \sTL_n$ is the disjoint union
of the sets $S(D)$ for $D \in \D_\max(I)$.  This follows from
Lemmas~\ref{lem:include} and \ref{lem:Imax}.
\end{proof}

Suppose $D \in \sTL_n$ is a (possibly odd) symmetric TL-diagram on
$4n$ vertices.  We define a subset $I(D) \subset [2n]$ by
$$
I(D) = \{ i \in [2n] \mid (i < j) \in D \; \text{or} \; (i,i') \in
D\}.
$$
Note that $|I(D)| = 2n - |D|$, so that $I(D)$ has even cardinality
whenever $D \in \sTL^e_n$.  Recall from before
Theorem~\ref{thm:stan} the definition of a standard partition of
$[2n]$.

\begin{lem}\label{lem:bij}
The map $D \mapsto (I(D),\overline{I(D)})$ is a bijection with image
equal to the set of standard partitions of $[2n]$ with at most $2$
parts.
\end{lem}
\begin{proof}
We describe how to recover $D$ from $I(D)$.  Let $\overline{I(D)} =
\{ j_1 < j_2 < \cdots < j_k\}$.  Then it must be the case that $(j_1
- 1,j_1) \in D$.  More generally suppose we know all the edges of
$D$ connected to $\{j_1,j_2, \ldots, j_{l-1}\}$ for some $l \leq k$.
Then $(i,j_l)$ is an edge of $D$, where $i \in I(D)$ is the maximum
number in $I(D)$ which is less than $j_l$ and which is not connected
to $\{j_1,j_2, \ldots, j_{l-1}\}$.  Furthermore, it is clear that
this algorithmic definition of the inverse map $(I,\bar I) \mapsto
D$ terminates successfully if and only if $(I, \bar I)$ is a
standard partitioning.

\end{proof}

\begin{corollary}
The dimension of $\CP_n$ is $2n-1 \choose n$.
\end{corollary}

\begin{proof}
This is an immediate corollary of Theorem \ref{thm:stan},
Proposition \ref{prop:size} and Lemma \ref{lem:bij}.
\end{proof}

Let $I, J \subset [2n]$ be two subsets of the same cardinality.  We
say $I = \{i_1 < \cdots< i_k\}$ is lexicographically smaller than $J
= \{j_1 < \cdots< j_k\}$ and write $I \prec_{\lex} J$ if for some $1
\leq l \leq k$ we have $i_1 = j_1, i_2 = j_2, \ldots, i_{l-1} =
j_{l-1}, i_l < j_l$.  We now define a total order $\prec$ on subsets
of $[2n]$. Suppose $I, J \subset [2n]$.  We define $I \prec J$ if
$|I| > |J|$ or $|I| = |J|$ and $I \prec_{\lex} J$.  We use the map
$D \mapsto I(D)$ to give an induced total order on $\sTL_n$: we have
$D \prec D'$ if $I(D) \prec I(D')$.

\begin{figure}[ht]
\begin{center}
\input{pfaf6.pstex_t}
\end{center}
\caption{The order $\prec_{\lex}$ on $\sTL_2$.} \label{pfaf6}
\end{figure}

\begin{lem}
\label{lem:comp} Let $D, D' \in \sTL_n$.  If $D$ is $I(D')$
compatible then $D \prec D'$.
\end{lem}
\begin{proof}
Suppose $D$ is $I(D')$ compatible.  Then $\overline{I(D')}$ must
have at least $|D| = |\overline{I(D)}|$ elements, so we have
$|I(D')| \leq |I(D)|$.  Thus we may suppose $|D| = |D'| = k$. Let
$\{(i_1 < j_1),\cdots, (i_k <j_k)\}$ be the vertical edges of $D$
(on the left side) and suppose that $j_1 < j_2 < \cdots < j_k$. If
$D$ is $I(D')$-compatible then $|\overline{I(D')} \cap (i_l,j_l)| =
1$ for each $l \in [1,k]$, so we must have $\overline{I(D')}
\prec_{\lex} \overline{I(D)}$.  This in turn implies that $I(D)
\prec_{\lex} I(D')$, so $D \prec D'$.
\end{proof}

\begin{example}
For $n=2$ we get $\{1,2,3,4\} \prec_{\lex} \{1,2,3\} \prec_{\lex}
\{1,2,4\} \prec_{\lex} \{1,3,4\} \prec_{\lex} \{1,2\} \prec_{\lex}
\{1,3\}$ which gives us the order on $\sTL_2$ as shown in
Figure~\ref{pfaf6}.
\end{example}

\begin{proposition}
\label{prop:tri} The transition matrix (given by
Theorem~\ref{thm:decomp}) from the set $\{\pf_{I, \bar I} \mid (I,
\bar I) \; \text{is standard}\}$ of standard complementary pfaffians
to the set $\{\Pfaf_D(A) \mid D \in \sTL^e_n\}$ is upper triangular
with 1's on the diagonal, under the order $\prec$.
\end{proposition}


\begin{proof}
Clearly $D \in \D_\max(I(D))$ so the matrix of the Proposition has
1's along the diagonal.  Since $\D_\max(I) \subset \D(I)$, by
Theorem~\ref{thm:diag} the coefficient of $\Pfaf_D$ in $\pf_{I, \bar
I}$ is non-zero if and only if $D$ is $I$-compatible.  By
Lemma~\ref{lem:comp}, $D$ is $I(D')$ compatible only if $D \prec
D'$, giving the upper triangularity.
\end{proof}

\begin{example}
For $n=2$ one obtains the transition matrix
$$
\left(
\begin{matrix}
1 & 1 & 1\\
0 & 1 & 1\\
0 & 0 & 1
\end{matrix}
\right)
.$$

We have labeled the rows by the standard complementary pfaffians
$(I, \bar I) = (\{1,3\},\{2,4\})$, $(\{1,2\},\{3,4\})$ and
$(\{1,2,3,4\},\{\emptyset\})$ from top to bottom and we label the
columns by the symmetric even TL-diagrams with vertical edges
$\{(1,2), (3,4)\}$, $\{(1,4), (2,3)\}$, $\{\emptyset\}$ from left to
right.
\end{example}

\begin{theorem}
\label{thm:basis} The TL-pfaffinants $\{\Pfaf_D(A) \mid D \in
\sTL^e_n\}$ form a basis for $\CP_n$.
\end{theorem}
\begin{proof}
This follows from Theorem~\ref{thm:stan} and
Proposition~\ref{prop:tri}.
\end{proof}
We will obtain another proof of Theorem~\ref{thm:basis} in
Section~\ref{sec:linind}.

\begin{problem}
Do the diagram pfaffinants $\{\Pfaf'_D(A) \mid D \in \sTL_n\}$
always lie in $\CP_n$?  If so, how are they expressed in the basis
of $TL$-pfaffinants and in the basis of standard complementary
pfaffians?
\end{problem}

By Examples \ref{ex:diag} and \ref{ex:even} the answer to the first
question is affirmative for $n=2$.  Note also that by
Proposition~\ref{prop:size} the number of diagram pfaffinants is
twice larger than the dimension of $\CP_n$, so if the diagram
pfaffinants $\{\Pfaf'_D(A)\}$ do lie in $\CP_n$ there must be
non-trivial relations among them.

\section{Pfaffians and non-intersecting paths in networks} \label{netw}
\subsection{Stembridge's network interpretation of Pfaffians}
John Stembridge in \cite{St} introduced an interpretation of
pfaffians in terms of networks. Let $G = (V, E)$ be a finite acyclic
directed graph.  We say that two directed paths in $G$ {\it
{intersect}} if they have a common vertex.  If $W$ and $U$ are
ordered sets of vertices of $G$, we say that $W$ is {\it
$G$-compatible} with $U$ if whenever $u<u'$ in $W$ and $v>v'$ in
$U$, every path from $u$ to $v$ intersects every path from $u'$ to
$v'$.


Let us suppose that a {\it {weight function}} $w: E \longrightarrow R$,
where $R$ is some ring, has been fixed.  For a $G$-path $p$, let
$w(p) = \prod_{e \in p} w(e)$ where the product is taken over all
edges in $p$.  For $u \in V$, $W \subset V$ let $P(u, W)$ denote the
set of $G$-paths from $u$ to any $v \in I$, and let $Q(u,W)$ be the
associated weight function $Q(u,W) = \sum_{p \in P(u, W)} w(p)$.
Similarly, for an $r$-tuple ${\bf u} = (u_1, \ldots, u_r)$ let
$P({\bf u}, W)$ denote the set of $r$-tuples of paths $(p_1, \ldots,
p_r)$ such that $p_i \in P(u_i, W)$.   The weight $w(p_1, \ldots,
p_r)$ of a $r$-tuple of paths is the product of the weights of each of the paths.  Let $P_0(\u, W) \subset P({\bf u}, W)$ denote the subset of non-intersecting tuples of paths.  We define $Q(\u,W) = Q_0(\u, W)$ to be  the sum of the weights of the elements of $P_0({\bf u}, W)$.

\begin{theorem}[{\cite[Theorem 3.1]{St}}]
\label{thm:Ste}
Let ${\bf u} = (u_1, \ldots, u_r)$ be an $r$-tuple of vertices in an
acyclic digraph $G$, and assume that $r$ is even. If $W \subset V$ is an ordered
subset of vertices such that $\bf u$ is $G$-compatible with $W$,
then
$$Q(\u,W) = \pf\left([Q((u_i,u_j),W)]_{1 \leq i < j \leq r}\right).$$
\end{theorem}

For convenience, if $G$ is an acyclic directed graph and ordered
vertex sets $\u = (u_1,\ldots,u_{2n}) \subset V$ and $W \subset V$
have been chosen we call the triple $N = (G,\u,I)$ a {\it network}.
For a network $N$, we define $P(N) = P(\u,W)$ and $P_0(N) =
P_0(\u,W)$.  We also let $Q(N)$ denote the weight sum $Q(\u,W)$, and
let $A(N) = A(G,\u,I)$ denote the array $(a_{ij} =
Q(\{u_i,u_j\},W)_{1 \leq i < j \leq r})$.

If $I \subset [2n]$, we let $u_I = \{u_i\}_{i \in I}$ denote the
corresponding set of vertices.  We then set $P_I(N) \subset P(N)$ to
be the subset of paths $\p = (p_1,\ldots,p_{2n})$ such that $p_i$
and $p_j$ do not intersect if both $i,j \in I$ or both $i,j \in \bar
I$.  We call the paths $\p \in P_I(N)$ compatible with $I$.  Thus
$P_0(N) = P_\emptyset(N) = P_{[2n]}(N)$.  We finally define $Q_I(N)$
to be the sum of the weights of the paths in $P_I(N)$.

The following statement is immediate from Theorem~\ref{thm:Ste} and the definitions we have made.
\begin{corollary}
\label{cor:Ste} Let $N = (G, \u, W)$ be a $G$-compatible network and
$I \subset [2n]$ be of even cardinality.  Then
$$
Q_I(N) = \pf_{I, \bar I}(A(N)).
$$

\end{corollary}

\subsection{Planar network definition of Pfaffinants}
Let $N = (G, \u, W)$ be a fixed network.  We assume that $G$ is planar and that a Jordan curve
$C$ passes through the sets $\u$ and $W$ of vertices so that $G$ is contained completely in the
interior of $C$.  We also assume that $\u$ and $W$ are contained in disjoint segments of $C$ so
that the ordering of $\u$ and $W$ is consistent with the arrangement of these vertices on $C$.  With this assumption, the $G$-compatibility of $\u$ and $W$ is immediate.  For short we will call a network $N$ satisfying these assumptions a {\it planar network}.

Suppose that $\p= (p_1, p_2, \ldots, p_{2n}) \in P(N)$ is a family
of paths such that no three paths in $\p$ intersect at the same
vertex.    Removing all the edges of $N$ that do not lie on any of
the paths $p_i \in \p$, and in addition marking all the edges of $N$
used twice by $\p$ we obtain a {\it marked network} $\tN = \tN(\p)$.
Note that by our assumption an edge of $N$ can be used at most twice
by the path family $\p$.  We say that $\p$ {\it covers} $\tN$ and
denote the set of coverings of $\tN$ by $P(\tN)$.  If $\tN$ is the
marked network obtained from some $\p \in P(N)$ we call $\tN$ a
marked subnetwork of $N$ and write $\tN \ll N$.  The weight $w(\tN)$
of a marked subnetwork is the weight $w(\p)$ for any path family
covering $\tN$.

\begin{figure}[ht]
\begin{center}
\input{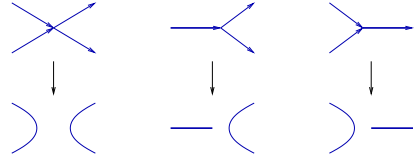}
\end{center}
\caption{The vertical uncrossing of an intersection.} \label{pfaf7}
\end{figure}

Suppose $p_i$ and $p_j$ intersect at some vertex $v$.  Then there
are two (possibly not distinct) edges $e_i \in p_i, e_j \in p_j$
entering $v$ and two edges $f_i \in p_i$ and $f_j \in p_j$ leaving
$v$. The {\it vertical uncrossing} of $v$ is obtained by detaching
$v$ into two new vertices $v_e$ and $v_f$ so that $v_e$ is incident
with $e_i$ and $e_j$ while $v_f$ is incident with $f_i$ and $f_j$,
as it is illustrated on Figure \ref{pfaf7}.  Alternatively, if the
vertices $\u$ are arranged on the left, the vertices $W$ arranged on
the right, and all edges are directed strictly from left to right,
then the vertical uncrossings always look like ``$\noXv$''.

Define an undirected graph $\Theta(\tN)$ by vertically uncrossing
every intersection point of $\tN$, removing all the marked edges and
ignoring all the orientations.  Note that $\Theta(\tN)$ does not
depend on $\p$, only on $\tN$.  The graph $\Theta(\tN)$ is a
disjoint union of a number of cycles, together with a number of
paths. We define the multiplicity of the marked network $\tN$ by
$\mult(\tN) = 2^{r}$ where $r$ is equal to the number of connected
components of $\Theta(\tN)$ which do not contain any of the vertices
in $\u$.

The components of $\Theta(\tN)$ containing one or more of the
vertices of $\u$ are a collection of paths which give rise to a
matching $\type(\tN)$ of $[2n] \cup [2n]'$: if $u_i,u_j$ belong to
the same component of $\Theta(\tN)$ then $(i,j), (i',j') \in
\type(\tN)$.  If $u_i$ does not belong in any component with some
other $u_j$, then $(i,i') \in \type(\tN)$.

\begin{lem}
\label{lem:Jordan}
Let $\p \in P(N)$ be a family of paths such that no three paths in $\p$ intersect at the same vertex and let $\tN = \tN(\p)$.  Then
$\type(\tN) \in \sTL_n$.
\end{lem}
\begin{proof}
We need to check that if $(i,j) \in \type(\tN)$ and $i < k < j$ then
$(k,l) \in \type(\tN)$ for some $i < l < j$.  The components of
$\Theta(\tN)$ are simple curves in the interior of the Jordan curve
$C$ connecting two points on the boundary of $C$.  The assumption
that $\u$ is arranged in order along the boundary of $C$ immediately
implies the required criterion.
\end{proof}

The definition of $\Theta(\tN)$ does not rely on the assumption that
the graph is drawn inside a Jordan curve, but Lemma~\ref{lem:Jordan}
does.

\begin{lem}
\label{lem:mult}
Let $N$ be a planar network and $\tN \ll N$ a marked subnetwork.  Suppose $I \subset [2n]$.  Then
the number of
path families which cover $\tN$ and are $I$-compatible  is given by
$$
|P(\tN) \cap P_I(N)| = \begin{cases} \mult(\tN) & \mbox{if $\type(\tN) \in \D(I)$,}
                            \\ 0 & \mbox{otherwise.} \end{cases}
$$
In particular, $|P(\tN) \cap P_I(N)|$ only depends on whether there is some $I$-compatible path family covering $\tN$.
\end{lem}
\begin{proof}
For each $\p \in P(\tN)$ we orient $\Theta(\tN)$ in the following manner.  If an edge $e \in \Theta(\tN)$ belongs to $p_i$ where $i \in I$
we orient $e$ with the same orientation as in $N$, that is, from $\u$ to $W$.   If an edge $e \in \Theta(\tN)$ belongs to $p_j$ where $j \in \bar I$
we orient $e$ with the opposite orientation to the one in $N$.  Since we removed all the marked edges when we produced $\Theta(\tN)$ no edge $e \in \Theta(\tN)$ receives both orientations.

The resulting directed graph $\Theta(\tN)_\p$ is a disjoint union of
directed paths and directed cycles.  This follows from the fact that
every intersection of $\tN$ involves a pair of paths $(p_i,p_j)$
where $i \in I$ and $j \in \bar I$.  One now checks that  $\p
\mapsto \Theta(\tN)_\p$ is a bijection between path families in $\p
\in P(\tN)$  and such directed graphs.

In addition, $\p \in P_I(N)$ if and only if the directed path in $\Theta(\tN)_\p$ that $u_i$ lies on is directed away from $u_i$ if $i \in I$ and directed towards $u_i$ if $i \in \bar I$.  This requirement can be satisfied only if $\type(\tN) \in \D(I)$.  The number of orientations of $\Theta(\tN)$ satisfying this additional condition is by definition equal to $\mult(\tN)$.
\end{proof}

For $D \in \sTL_n$ define the following function $\hat{\Pfaf}'_D: \{\text{planar networks}\} \to R$ on planar networks:
$$\hat{\Pfaf}'_{D}(N) = \sum_{\substack{\tN \ll N\\ {\type}(\tN)=D}}
\mult(\tN)w(\tN).$$

\begin{proposition}
\label{prop:net}
Let $I \subset [2n]$ be of even cardinality and $N$ be a planar network.  Then
$$
\pf_{I,\bar I}(A(N)) = \sum_{D \in \D(I)} \hat{\Pfaf}'_D(N).
$$
\end{proposition}

\begin{proof}
By Corollary~\ref{cor:Ste}, $\pf_{I,\bar I}(A(N))$ is the sum of the weights of the $I$-compatible families of paths $P_I(N)$.  Thus
\begin{align*}
\pf_{I, \bar I}(A(N)) & = \sum_{\p \in P_I(N)} w(\p) \\
&= \sum_{\tN \ll N} \left(\sum_{\p \in P_I(N) \cap P(\tN)} w(\p) \right) \\
& = \sum_{\substack{\tN \ll N \\ \type(\tN) \in \D(I)}} \mult(\tN)w(\tN) & \mbox{by Lemma~\ref{lem:mult}} \\
&=  \sum_{D \in \D(I)} \hat{\Pfaf}'_D(N).
\end{align*}



\end{proof}

Now for $D \in \sTL^e_n$ define $\hat{\Pfaf}_D: \{\text{planar networks}\} \to R$ by
$$\hat{\Pfaf}_{D}(N) = \sum_{D' \in S(D)}
\hat{\Pfaf}'_{D'}(N).$$

\begin{theorem}
\label{thm:equal}
Let $D \in \sTL^e_n$ and $N$ be a planar network.  Then
$$
\Pfaf_D(A(N)) = \hat{\Pfaf}_D(N).$$
\end{theorem}
\begin{proof}
The proof of Theorem~\ref{thm:decomp} and Proposition~\ref{prop:net} shows that
$$\pf_{I, \bar I}(A(N)) = \sum_{D \in \D_\max(I)} \hat{\Pfaf}_D(N).$$  Using the statement and the proof of Proposition~\ref{prop:tri}
we see that $\hat{\Pfaf}_D(N)$ and $\Pfaf_D(A(N))$ can be expressed in terms of $\pf_{I,\bar I}(A(N))$ in an identical way so we
conclude that $\Pfaf_D(A(N)) = \hat{\Pfaf}_D(N)$.
\end{proof}

\begin{remark}
Observe that functions $\hat \Pfaf'_D(N)$ do {\it {not}} coincide
with the evaluations $\Pfaf'_D(A(N))$ of diagram pfaffinants.  In
particular the diagram pfaffinants $\Pfaf'_D(A)$ might take negative
values when evaluated at $A(N)$ for a planar network $N$.
\end{remark}

\subsection{Independence of Temperley-Lieb pfaffinants}
\label{sec:linind} We will show directly using
Theorem~\ref{thm:equal} that the elements $\{\Pfaf_D(A) \mid D \in
\sTL^e_n\}$ are linearly independent.  This will give us alternative
proofs of Theorems~\ref{thm:stan} and \ref{thm:basis}.

Let $D \in \sTL_n$.  We will now define a planar network $N(D)$ with
the property that $\hat{\Pfaf'}_{D'}(N(D))$ is non-zero if and only
if $D = D'$.  The network $N(D)$ is embedded into the plane $\R^2$
in a particular way. First, place the vertices $u_1,\ldots, u_{2n}$
on the line $x = 0$ so that $u_i$ has coordinates $(0,2n-i)$.  For
an edge $(i<j) \in D$ we call the vertex $i$ {\it {outgoing}} and
the vertex $j$ {\it {ingoing}}. The vertices $i$ such that $(i,i')
\in D$ are neither outgoing nor ingoing. Now place the ``sink''
vertices $W$ as follows: for each $i \in [2n]$ such that $(i,i') \in
D$ or $(i < j) \in D$ we place $w_i \in W$ at coordinates
$(1,2n-i)$.  To obtain the rest of $N(D)$, we first join $u_i$ with
$w_i$ with a straight line whenever $w_i$ exists, that is when $i$
is not ingoing.  Finally we join $u_{j_k}$ with $w_{i_k}$ where $j_1
< j_2 < \cdots $ are the ingoing vertices and $i_1 < i_2 < \cdots$
are the outgoing vertices.  The intersection of any of these lines
is also defined to be a vertex of $N(D)$ which does not belong to
either $\u$ or to $W$.  All edges are directed so that the
$x$-coordinate increases along the edges.

\begin{figure}[ht]
\begin{center}
\input{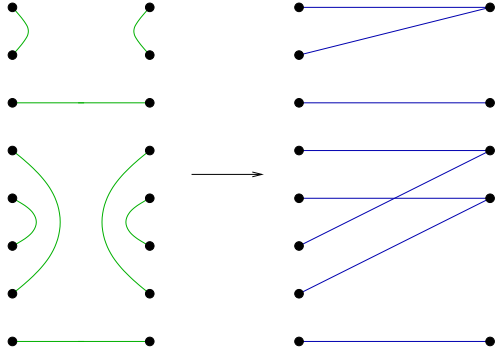}
\end{center}
\caption{A symmetric TL-diagram $D$ and the corresponding network
$N(D)$} \label{fig:pfaf4}
\end{figure}

Note that no three of the drawn lines intersect at one point, since
by construction the set of these lines is a union of two pairwise
non-intersecting families of lines.  An example of this construction
of $N(D)$ is shown in Figure~\ref{fig:pfaf4}.

Say that an edge $(i,j) \in D$ is on the {\it outside} if one cannot
find $(k,l) \in D$ so that $1 \leq k < i < j < l \leq 2n$.  The
network $N(D)$ is the union of the networks $N(D_{[i,j]})$ for
outside edges $(i,j)$ together with the networks $N(D_{i})$ for
horizontal edges $(i,i')$.  Here $D_{[i,j]}$ denotes the obvious
restriction of a diagram $D$ to the set of vertices $[i,j] \cup
[i',j'] \subset [2n] \cup [2n]'$. Let $\tN(D) \ll N(D)$ denote the
marked subnetwork consisting of all edges of $N(D)$, with no edges
marked.

\begin{lemma} \label{lem:constr}
We have $\type({\tN(D)}) = D$.  Let $\p \in P(N)$ be a family of
paths such that no three paths intersect at the same vertex. Then
$\tN(\p) = \tN(D)$.
\end{lemma}

\begin{proof}
By the previous comments, it is enough to prove the lemma for each
of the networks $N(D_{[i,j]})$ corresponding to outside edges $(i,j)
\in D$.  We proceed by induction on $|j-i|$, the base case being
trivial.  All vertices in $[i,j]$ are outgoing or ingoing, and there
are twice as many source vertices $\u$ as sink vertices $W$ in
$N(D_{[i,j]})$. Call the edges of $N(D_{[i,j]})$ incident to the
sink vertices the {\it {outer skeleton}} $Sk(N(D_{[i,j]}))$.

Now remove the outer skeleton from $N(D_{[i,j]})$.  We obtain a
network $N(D_{[i,j]})'$ isomorphic to $N(D_{[i+1,j-1]})$, which is
the union of the networks $N(D_{[i_p,j_p]})$, where $\{(i_p,j_p)\}$
is the set of outside edges formed when we remove edge $(i,j)$ from
$D_{[i,j]}$.  Under this identification, the sink vertices of
$N(D_{[i,j]})'$ are the intersection points of the pairs of segments
$\{(u_{j_k}, w_{i_k}),(u_{i_{k+1}}, w_{i_{k+1}})\}$.  By the
inductive assumption, we have $\type(\tN(D_{[i+1,j-1]})) =
D_{[i+1,j-1]}$ and since $Sk(N(D_{[i,j]}))$ (after redirecting the
edges) is a path from $u_i$ to $u_j$, it follows immediately that
$\type({\tN(D_{[i,j]})}) = D_{[i,j]}$.

By the inductive assumption applied to each $N(D_{[i_p,j_p]})$,
there is only one marked network of $N(D_{[i,j]})'$ arising from a
family of paths $\p \in P(N)$ without triple intersections.  Each of
the sink vertices of $N(D_{[i,j]})'$ has incoming degree 2, and thus
$\p$ must cover (counted with multiplicity) two of the outgoing
edges from each such vertex.  However, $\p$ must contain the two
paths consisting of the single edge $(u_i,w_i)$ and the single edge
$(u_{j_s},w_{i_s})$, where $j_s = j$.  A simple counting argument
shows that each sink vertex $w_{i_r}$ is incident with exactly two
paths.  Combining these facts, one concludes that each edge of
$Sk(N(D_{[i,j]}))$ is covered by $\p$ exactly once.
\end{proof}

\begin{figure}[ht]
\begin{center}
\input{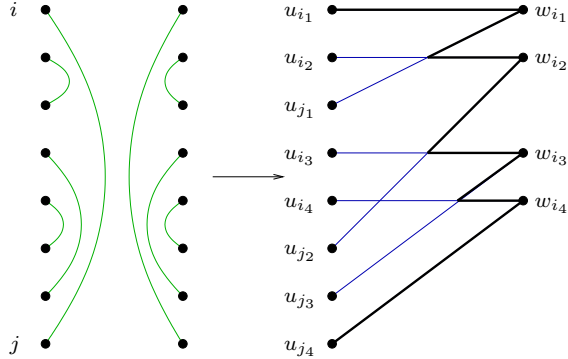}
\end{center}
\caption{A symmetric TL-diagram $D$, with corresponding network
$N(D)$ and outer skeleton $Sk(N(D))$ shown in bold.}
\label{fig:pfaf5}
\end{figure}

An illustration of the proof is shown in Figure \ref{fig:pfaf5}.

\begin{theorem} \label{thm:linind}
The elements $\{\Pfaf_D(A) \mid D \in \sTL^e_n\} \subset \CP_n$ are
linearly independent.
\end{theorem}

\begin{proof}
Suppose there is a non-trivial linear combination $c = \sum_{D \in
\sTL^e_n} c_D \Pfaf_D(A)$ of the $\Pfaf_D(A)$-s which evaluates to
$0$ for any upper triangular array $A$.  Then in particular it
should evaluate to $0$ on $A(N)$ for a planar network $N$.  Let $D
\in \sTL^e_n$ be such that $\Pfaf_D(A)$ enters the expression with a
non-zero coefficient $c_D$, and $|D|$ is the largest possible
satisfying this condition. Then by Lemma~\ref{lem:constr},
$\hat{\Pfaf'}_{D}(N(D))$ contributes a non-zero value to
$\Pfaf_D(A(N(D)))$ but we have $\hat{\Pfaf'}_{D'}(N(D)) = 0$ for all
other $D' \neq D$ such that $D' \in S(D)$.  However, by the choice
of $D$ we have $\Pfaf_{D'}(A(N(D))) = 0$ for all other $D'$ such
that $c_{D'} \neq 0$.  We conclude $c_D = 0$, obtaining a
contradiction.
\end{proof}

Theorem~\ref{thm:linind} gives alternative proofs of
Theorems~\ref{thm:stan} and \ref{thm:basis} without relying on
results of \cite{DP}.

\subsection{Network positivity}
\label{sec:npos} Call a skew symmetric matrix $A$ {\it
{network-positive}} if it is equal to $A(N)$ for some planar network
$N$ with positive weights on edges (we assume the coefficient ring
$R = {\mathbb R}$).

The notion of network positivity is a substitute for the notion of
{\it {total non-negativity}} of matrices.  Recall that an arbitrary
matrix $M$ is totally non-negative if all its minors are
non-negative.  It is known (see for example \cite[Theorem 3.1]{Br})
that every totally non-negative matrix arises from a planar network.

It is not clear how to make a similar definition for skew-symmetric
matrices.  The following example is taken from \cite{Kim}. Take the following skew-symmetric matrix:
$$A =
\left(
\begin{matrix}
0 & 1 & 0 & 0\\
-1 & 0 & 1 & 0\\
0 & -1 & 0 & 0\\
0 & 0 & 0 & 0
\end{matrix}
\right).
$$

Every skew-symmetric submatrix of $A$ of even size has a
non-negative pfaffian.  However, as we will now show, $A$ is not
equal to $A(N)$ for any planar network $N$.  Thus the naive
generalization does not seem to be appropriate.

\begin{lemma}\label{hyun}
$A$ is not equal to $A(N)$ for any positive planar network $N$.
\end{lemma}

\begin{proof}
Indeed, assume $\bf u$ and $W$ are placed on the boundary of a
Jordan curve.  Since $a_{23} \not = 0$ there should be a pair of
non-intersecting paths $p_2$ and $p_3$ from $u_2$ and $u_3$ to $W$
(see Figure \ref{fig:pfaf12}).

\begin{figure}[ht]
\begin{center}
\input{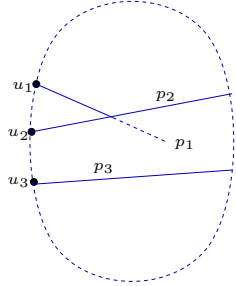}
\end{center}
\caption{It is impossible to have $a_{12}, a_{23}>0$ while
$a_{13}=0$ in $A(N)$ for a non-negative planar network $N$.}
\label{fig:pfaf12}
\end{figure}

Since $a_{12} \not = 0$ there should be at least one path $p_1$ from
$u_1$ to $W$. Since $a_{13}=0$, the path $p_1$ must intersect $p_3$,
and therefore $p_2$. However, in that case if we traverse $p_1$ up
to the point of intersection with $p_2$ and continue along $p_2$, we
obtain a path from $u_1$ to $W$ not intersecting $p_3$,
contradicting our assumptions.
\end{proof}

\subsection{Network positivity and pfaffinants}

\begin{proposition}
For a network-positive $A$ and any $D \in \sTL^e_n$ we have
$\Pfaf_{D}(A) \geq 0$.
\end{proposition}

\begin{proof}
We know from Theorem \ref{thm:equal} that $\Pfaf_{D}(A)$ has an
interpretation as the weight-multiplicity generating function of
certain marked subnetworks of $N$. The statement follows
immediately.
\end{proof}


For any $K \in \CP_n$ one can formally write $K$ as a linear
combination of the symbols $\Pfaf'_{D}$.  Namely, by
Theorem~\ref{thm:basis} one can express $K = \sum c_D \Pfaf_{D}$ in
terms of TL-pfaffinants. Now we use the expansions $\Pfaf_{D}(A) =
\sum_{D' \in S(D)} \Pfaf'_{D'}(A)$ to obtain the needed formal
presentation $K = \sum c'_D \Pfaf'_{D}$.

\begin{theorem} [cf. Corollary 3.6, \cite{RS1}] \label{thm:npos}
Let $K \in (\CP_n)_\R$. The following are equivalent:
\begin{enumerate}
 \item for any network-positive $A$ one has $K(A) \geq 0$;
 \item The coefficients $c'_D$ in $K = \sum_{D \in \sTL_n} c'_D \Pfaf'_{D}$ are non-negative.
\end{enumerate}
\end{theorem}
We call an element $f \in \CP_n$ {\it network positive} if it
satisfies one of the conditions (and thus both) of
Theorem~\ref{thm:npos}.

\begin{proof}
By Theorem~\ref{thm:equal}, $K(A(N)) = \sum c'_D \hat \Pfaf'_{D}(N)$
and each $\hat \Pfaf'_{D}(N)$ by definition enumerates sums of
weights of certain marked subnetworks of $N$, one direction is
obvious. It was shown in Lemma~\ref{lem:constr} that the networks $N(D)$
possess the property that $\hat \Pfaf'_{D'}(N(D)) \neq 0$ if and
only if $D'= D$. This implies the other direction.
\end{proof}

Let $C_n \subset \CP_n$ denote the cone consisting of network
positive elements.  Theorem~\ref{thm:npos} shows that $C_n$ is
rational and polyhedral and a simple argument using the networks
$N(D)$ shows that $C_n$ is pointed (contains no lines).  However,
the cone $C_n$ possesses some interesting polyhedral geometry and
the edge generators of $C_n$ are rather tricky to describe.  Finding
generators of the semigroup $C_n \cap \Z[\pf_{I, \bar I} \mid I \subset [2n]]$
of integral points is even trickier.  Note that by
Theorem~\ref{thm:decomp} and Proposition~\ref{prop:tri}, the
$\Z$-span of $\{\Pfaf_D(A) \mid D \in \sTL^e_n\}$ is equal to the
$\Z$-span of $\{\pf_{I, \bar I} \mid I \subset [2n]\}$.

The description of the edge generators of $C_n$ can be simplified to
a combinatorial problem concerning boolean lattices.

Let us call an even symmetric diagram $D \in \sTL^e_n$ {\it maximal}
if it is $I$-maximal for the subset $I = I_{\rm
alt}=\{1,3,5,\ldots,2n-1\}$.  Since $\D(I_{\rm alt}) = \sTL^e_n$, a
diagram $D \in \sTL^e_n$ is maximal if no odd edges can be added to
it.  By Lemma~\ref{lem:Imax}, $D \in \sTL^e_n$ is maximal if and
only if for every $D'$ so that $D \in S(D')$ we have $D = D'$.  The
following Lemma says that to find the edge generators of $C_n$ we
may restrict our attention to elements $f \in \CP_n$ which are
linear combinations of TL-pfaffinants labeled by a set $S(D_m) \cap
\sTL^e_n$ for maximal $D_m$.

\begin{lem} \label{lem:impos}
Suppose $f = \sum_{D \in \sTL^e_n} c_D \Pfaf_D(A) \in C_n$ lies in
the network positive cone.  Then so does $f_{D_m} = \sum_{D \in
S(D_m)} c_D \Pfaf_D(A)$ for each maximal diagram $D_m \in \sTL^e_n$.
\end{lem}
\begin{proof}
Suppose when expressed in terms of diagram pfaffinants as in
Theorem~\ref{thm:npos} we have $f = \sum_{D' \in \sTL_n} c_{D'}
\Pfaf'_D(A)$ where $c_{D'} \geq 0$ is given by
\begin{equation}
\label{eq:DprimeD} c_{D'} = \sum_{\substack{D \in \sTL^e_n \\D' \in
S(D)}} c_D.
\end{equation}
 Suppose $D_m$ is maximal and $D \in S(D_m)$.  By
Lemma~\ref{lem:Imax}, the summation in (\ref{eq:DprimeD}) can be
taken over $D \in (\sTL^e_n \cap S(D_m))$ satisfying $D' \in S(D)$
instead.  Also using Lemma~\ref{lem:include}, this shows that
$f_{D_m}$ lies in $C_n$.
\end{proof}

Now let $D_m$ be maximal.  By the proof of Lemma~\ref{lem:include}
the diagrams $D' \in S(D_m)$ form a boolean lattice $B_s = 2^{[s]}$
under the order $D_1 < D_2 \Leftrightarrow D_1 \in S(D_2)$.  When
$s$ is even, the even diagrams $S(D_m) \cap \sTL^e_n$ correspond to
the even levels $B_s^e$ in $B_s$.  When $s$ is odd, the even
diagrams $S(D_m) \cap \sTL^e_n$ correspond to the odd levels $B_s^o$
in $B_s$. The edge generators of $C_n$ can then be calculated by
solving the following problem.

\begin{problem}
Let $B_s$ be a boolean lattice and $\alpha \in \{o,e\}$.  What is
the cone of sequences $\{a_S\}_{S \in B_s^\alpha}$ of real numbers,
indexed by either the odd or the even subsets of $\{1,2,\ldots,
s\}$, which satisfies the condition
$$
\sum_{S' \subset S} a_S > 0
$$
for every $S' \in B_s$?
\end{problem}

\begin{example}\label{ex:boolean}
Assign to each element of the lattice $B_3 = 2^{\{a,b,c\}}$ one of
the formal variables $\emptyset, a, b, c, ab, ac, bc, abc$. Also,
define ${\bf a} = a + \emptyset$, ${\bf b} = b + \emptyset$, ${\bf
c} = c + \emptyset$, ${\bf abc} = abc + ab + ac + bc +a +b +c
+\emptyset$. We want to characterize the cone of $(t_1,t_2,t_3,t_4)
\in \R^4$ such that $t_{1} {\bf abc} + t_2 {\bf a} + t_3 {\bf b} +
t_4 {\bf c}$ is non-negative in terms of the original eight formal
variables.  It turns out that the edges of the this cone are
generated by the set $V_3$ of vectors $(0,1,0,0)$, $(0,0,1,0)$,
$(0,0,0,1)$, $(1,-1,-1,1)$, $(1,-1,1,-1)$, $(1,1,-1,-1)$.  However,
if we were to consider the problem restricted to $(t_1,t_2,t_3,t_4)
\in \Z^4$, we need to also add the vectors $(1,-1,0,0)$,
$(1,0,-1,0)$ and $(1,0,0,-1)$ to the above set.

Let $D \in \sTL^e_4$ be the diagram with vertical edges $(1,2),
(3,4), (5,8), (6,7)$. Label the (removable) odd edges $(1,2)$,
$(3,4)$ and $(5,8)$ with $a$, $b$ and $c$ correspondingly. Then
elements of $S(D)$ are in bijection with nodes of boolean algebra
$B_3$. Thus, a linear combination $\sum_{D' \in S(D)} t_{D'}
\Pfaf'_{D'}$ can be network-nonnegative if and only if the
coefficients of the four TL-pfaffinants corresponding to the nodes
$abc$, $a$, $b$ and $c$ of $B_3$ lies in the cone generated by the
set $V_3$.

\begin{figure}[ht]
\begin{center}
\input{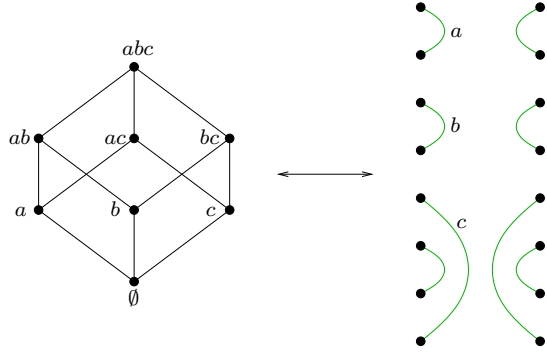}
\end{center}
\caption{The boolean algebra $B_3 = 2^{\{a,b,c\}}$ in
Example~\ref{ex:boolean}.} \label{fig:pfaf8}
\end{figure}

\end{example}

\subsection{Positive differences of complementary pfaffians}
\label{sec:comppfaf} One obtains a criterion for a linear
combination of complementary pfaffians $\pf_{I,\bar I}$ to be
network positive by combining Theorems~\ref{thm:decomp} and
\ref{thm:npos}.  The next result gives one way to produce network
positive differences of two complementary pfaffians.

First suppose $I \subset [2n]$ is an even subset and suppose that
$|I| \geq n$.  Suppose $I = \{i_1 < i_2 < \cdots < i_k\}$ and $\bar
I = \{j_1 < j_2 < \cdots < j_r\}$ where $r \leq k$.  We set
$\min(I,\bar I) =
\{\min(i_1,j_1),\min(i_2,j_2),\ldots,\min(i_r,j_r),i_{r+1},\ldots,i_k\}$.
For convenience, we may let $j_{r+1} = \cdots = j_{k} = \infty$.

\begin{proposition}
\label{prop:comppfaf} The difference $\pf_{\min(I,\bar
I),\overline{\min(I,\bar I)}} - \pf_{I, \bar I}$ is network
positive.
\end{proposition}
\begin{proof}
We shall show that $\D(I) \subset \D(\min(I,\bar I))$.  The result
will then follow from Theorems~\ref{thm:diag} and~\ref{thm:npos}. So
let $D \in \D(I)$ and suppose that $(i < j) \in D$.  Then either $i
\in I$ and $j \in \bar I$ or $i \in \bar I$ and $j \in I$.  We need
to show that exactly one of $(i,j)$ lies in $\min(I,\bar I)$.  The
key fact is that \begin{equation} \label{eq:same} |[i+1,j-1] \cap I|
= |[i+1,j-1] \cap \bar I|. \end{equation} Suppose that $i = i_a \in
I$ and $j = j_b \in \bar I$.  If $i_a < j_a$ then $i \in \min(I,
\bar I)$ and furthermore $i_b < j_b$ by (\ref{eq:same}) so that $j
\notin \min(I,\bar I)$.  Otherwise if $i_a > j_a$ we deduce by
(\ref{eq:same}) that $i_b > j_b$; so we conclude again that exactly
one of $(i,j)$ lies in $\min(I,\bar I)$.  The case that $i \in \bar
I$ and $j \in I$ is similar.
\end{proof}

%
%


\section{Relation between pfaffinants and immanants} \label{rel}
\subsection{Rhoades and Skandera's Temperley-Lieb immanants}
The {\it {Temperley-Lieb immanants}} were discovered by Rhoades and
Skandera~\cite{RS1}, who gave a number of remarkable positivity
properties of these immanants.  The exposition we now give is
similar to the presentation in~\cite{LPP} to which we refer for
unexplained notations.


Let $\TL_n$ be the set of Temperley-Lieb diagrams on $2n$ points
$\{1,2,\ldots,2n\}$, with $\{1,2,\ldots,n\}$ arranged top to bottom
on the left side of a rectangle and $\{n+1,\ldots,2n\}$ arranged
bottom to top on the right side. Let $w$ be a permutation in $S_n$.
By abuse of notation we also denote by $w$ a chosen wiring diagram,
thought of as a planar network connecting the $n$ source points on
the left to $n$ sink points on the right. Now uncross the crossings
of $w$ in all possible ways, each crossing becoming either a
vertical uncrossing ``\noXv'' or a horizontal uncrossing ``\noXh''.
Let $X(w)$ be the set of such uncrossings, and for $x \in X(w)$ let
$D(x)$ be the element of $\TL_n$ topologically equivalent to $x$
(with any loops removed). Let $h(x)$ be the number of horizontal
uncrossings in $x$ and let $l(x)$ be the number of loops formed.
Define the weight $\wt(x)$ of $x$ by $\wt(x) = 2^{l(x)}
(-1)^{h(x)}$.  For $d \in \TL_n$ define $f_{d}: S_n \to \Z$ by
$$f_d(w) = \sum_{\substack{x \in X(w)
\\ D(x) = d}} \wt(x).$$

Let $B = (b_{ij})$ be a $n \times n$ matrix.  Then for $d \in \TL_n$
the {\it {TL-immanant}} $\Imm_d(B)$ is defined as
$$\Imm_d(B) = \sum_{w \in S_n} f_d(w) b_{1, w(1)} \cdots b_{n,
w(n)}.$$

Let $S \subseteq [2n]$ and recall that $\bar S = [2n] \backslash S$
denotes the complement of $S$ in $[2n]$.  The {\it $S$-coloring} of
$[2n]$ is obtained by coloring the elements of $S$ black, and the
elements $\bar S$ white. We call a diagram $d \in \TL_n$ compatible
with $S$ (or simply $S$-compatible) if each edge of $d$ has ends of
different color in the $S$-coloring. We denote by $\D(S) \subset
\TL_n$ the set of $S$-compatible diagrams.

For two subsets $I, J\subset [n]$ of the same cardinality let $\Delta_{I,J}(B)$
denote the {\it minor\/} of an $n\times n$ matrix $B$ in the row set $I$ and the
column set $J$. Let $\hat I := [n]\setminus I$ and let $I^\wedge := \{2n+1-i\mid i\in I\}$.

\begin{theorem}
\label{thm:immdecomp} {\rm
Rhoades-Skandera~\cite[Proposition~4.3]{RS1}} \ For two subsets
$I,J\subset [n]$ of the same cardinality and $S=J\cup (\hat
I)^\wedge$, we have
$$
\Delta_{I,J}(B)\cdot \Delta_{\bar I, \bar J}(B)
= \sum_{d \in \D(S)} \Imm_d(B).
$$
\end{theorem}

\subsection{Expressing TL-immanants as TL-pfaffinants}

Let $A = (a_{ij})_{1 \leq i < j \leq 2n}$ be an uppertriangular
array such that $a_{ij} = 0$ if $1 \leq i < j \leq n$ or $n+1 \leq i
< j \leq 2n$.  Let $B = (b_{ij})$ be the $n \times n$ matrix given
by $b_{ij} = a_{i, j+n}$.  Our aim is to relate the TL-pfaffinants
$\Pfaf_D(A)$ with the TL-immanants $\Imm_d(B)$.

Call a set $I \subset [2n]$ balanced if $|I \cap [n]| =
\frac{|I|}{2}$.  Let $I_1 = I \cap [n]$ and $I_2 = \{i -n \mid i \in
I \cap [n+1,2n]\}$.

\begin{lemma}\label{lem:pfafdet}
Let $I \subset [2n]$ be an even subset.  Then $\pf_{I,\bar I}(A) =
0$ if $I$ (equivalently, $\bar I$) is not balanced. If $I$ is
balanced then $\pf_{I,\bar I }(A) = (-1)^{{|I_1| \choose 2}+{|\bar
I_1| \choose 2}} \Delta_{I_1,I_2}(B)\cdot \Delta_{\bar I_1, \bar
I_2}(B)$.
\end{lemma}

\begin{proof}
The first statement is clear since if $I$ is not balanced any
matching contains an edge corresponding to a zero entry of $A$. The
second statement follows from the observation that $\pf(A) = (-1)^{n
\choose 2} \Delta(B)$.

\end{proof}

Thus non-zero products of complementary pfaffians of $A$ are up to
sign equal to products of complementary minors of $B$.  Hence one
should be able to express the TL-pfaffinants of $A$ in terms of the
TL-immanants of $B$.

\begin{figure}[ht]
\begin{center}
\input{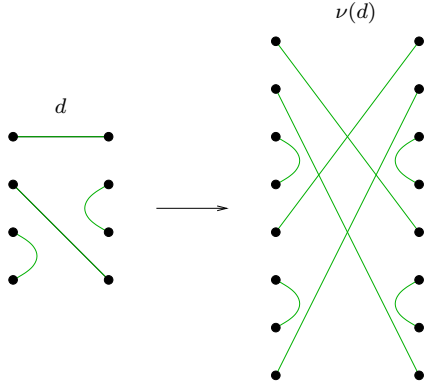}
\end{center}
\caption{The matching $\nu(d)$ produced from a TL-diagram $d$.}
\label{fig:pfaf9}
\end{figure}

Let $d \in \TL_n$.  Define a matching $\nu(d)$ of $[2n] \cup [2n]'$
as follows: interpret the left side of $d$ (originally labeled
$\{1,2,\ldots n\}$) as the vertices from $1$ to $n$ and the right
side of $d$ (originally labeled $\{2n,2n-1,\ldots,n+1\}$) as the
vertices from $(n+1)'$ to $2n'$.  Now force $\nu(d)$ to be
mirror-symmetric by adding the edge $(i,j)$ (resp. $(i',j')$,
$(i,j')$) whenever the edge $(i',j')$ (resp. $(i,j)$, $(i',j)$) is
present in $d$.  Let $X(d)$ be the set of all ways to uncross all
crossings in $\nu(d)$, where as in Section~\ref{sect:dpfaf} we
always uncross mirror symmetric crossings in the same manner.  As
usual, we pick the embedding of $\nu(\pi)$ so that no pair of edges
intersect more than once or have a point of tangency, and no three
edges intersect at a single point.

We define the {\it {weight}} $\wt(x)$ of an element $x \in X(d)$ as
$$
\wt(x) = 2^{l(x)} (-1)^{\uv(x)+ \ph(x)}
$$
where $l(x), \uv(x), \ph(x)$ are as defined in
Section~\ref{sect:dpfaf}.  Similarly we define $D(x) \in \sTL_n$ to
be the symmetric TL-diagram obtained from the uncrossing $x$.

We define $g_{D}: \TL_n \to \Z$ by
$$
g_D(d) = \sum_{\substack{x \in X(d) \\ D(x) = D}} \wt(x).
$$

Denote by $z(d)$ the number of edges in $d$ with both ends in $[n]$.
Finally, let $\tilde g_D(d) = (-1)^{z(d) \cdot n}g_D(d)$.

\begin{theorem} \label{thm:dec}
Let $D \in \sTL^e_n$.  Then
$$\Pfaf_{D}(A) = \sum_{D' \in S(D)}
\sum_{d \in \TL_n} \tilde g_{D'}(d) \Imm_d(B).$$
\end{theorem}

\begin{proof}
Let $\tilde \Pfaf_{D}(A)$ denote the right hand side of the equation
in the theorem:
$$\tilde \Pfaf_{D}(A) = \sum_{D' \in S(D)} \sum_{d \in \TL_n}
\tilde g_{D'}(d) \Imm_d(B).$$  By Proposition \ref{prop:tri}, it is
enough to show that the elements $\{\tilde \Pfaf_{D}(A)\}$ satisfy
the following decomposition formula (see Theorem~\ref{thm:decomp}):
$$\pf_{I, \bar I}(A) = \sum_{D \in \D_\max(I)} \tilde \Pfaf_D(A).$$

We have by definition
\begin{align*}
\sum_{D \in \D_\max(I)} \tilde \Pfaf_D(A) &= \sum_{D' \in \D(I)}
\sum_{d \in \TL_n} \tilde g_{D'}(d) \Imm_d(B) \\
&= \sum_{d \in \TL_n} (-1)^{z(d) \cdot n} \left(  \sum_{\substack{x
\in X(d) \\ D(x) \in \D(I)}} \wt(x) \right) \Imm_d(B).
\end{align*}
Now we proceed as in the proof of Theorem~\ref{thm:diag}.  Suppose
$x \in X(d)$ is an uncrossing of $\nu(d)$ such that $D(x) \in
\D(I)$.  We direct all the strands and loops in $x$ so that the
initial vertex of each strand belongs to $I \cup (\bar I)'$ (and,
thus the end vertex belongs to $\bar I \cup I'$). We allow the
closed loops to be directed in either direction.  Now define an
almost sign-reversing involution on this set of oriented diagrams
exactly as in Theorem~\ref{thm:diag}.

Thus the contribution of $\Imm_d(B)$ to $\sum_{D \in \D_\max(I)}
\tilde \Pfaf_D(A)$ is equal to the sum over the aligned uncrossings
$x \in X(d)$ of $(-1)^{z(d) \cdot n}(-1)^{\uv(x)+ \ph(x)}$.  As in
the proof of Theorem~\ref{thm:diag}, such an aligned uncrossing $x_d
\in X(d)$ is unique if it exists -- it corresponds to an orientation
$\mu(d)$ of $\nu(d)$ which connects elements of $I \cup (\bar I)'$
with elements of $\bar I \cup I'$.  Restricting $\mu(d)$ to
$\{1,2,\ldots,n\} \cup \{(n+1)', \ldots, (2n)'\}$, we see that $d$
must be $S = I_1 \cup (\hat{I_2})^\vee$-compatible, and in
particular $I$ must be balanced.  Conversely, if $d$ is $I_1 \cup
(\hat{I_2})^\vee$-compatible one obtains a unique such orientation
$\mu(d)$.

Finally we must calculate $(-1)^{\uv(x)+ \ph(x)}$ for $x(d)$.  The
unpaired crossings between $(i,j')$ and $(j,i')$ are always
uncrossed horizontally, so contribute nothing to the sign.  The
paired crossings which are uncrossed horizontally correspond to
pairs of edges $(i_1 < j_1) \in d$ and $(i_2 < j_2) \in d$, both of
which are horizontal and such that both $i_1,i_2 \in I_1$ or both
$i_1,i_2 \in I_2$.  Thus for $d \in \D(S)$ the coefficient of
$\Imm_d(B)$ in $\sum_{D \in \D_\max(I)} \tilde \Pfaf_D(A)$ is equal
to

$$(-1)^{z(d) \cdot n} (-1)^{{|I_1|-z(d) \choose 2}+{n - |I_1|-z(d) \choose
2}} = (-1)^{{|I_1|\choose 2}+{n - |I_1|\choose 2}}.$$ This identity
can be proven by induction on $z(d)$, noting that $(-1)^{k \choose
2} = 1$ if $k \equiv 0,1 \mod 4$ and $(-1)^{k \choose 2} = 1$ if $k
\equiv 2,3 \mod 4$.  Now summing over over all $d \in \TL_n$ and
using Theorem~\ref{thm:immdecomp} and Lemma~\ref{lem:pfafdet}, we
see that $\pf_{I, \bar I}(A) = \sum_{D \in \D_\max(I)} \tilde
\Pfaf_D(A)$.

\end{proof}

\begin{example}
Let us take a $4 \times 4$ skew-symmetric matrix
$$A =
\left(
\begin{matrix}
0 & 0 & x & y\\
0 & 0 & z & t\\
-x & -z & 0 & 0\\
-y & -t & 0 & 0
\end{matrix}
\right)
$$
and let $B$ be its $2 \times 2$ minor
$$B =
\left(
\begin{matrix}
x & y\\
z & t
\end{matrix}
\right)
.$$

\begin{figure}[ht]
\begin{center}
\input{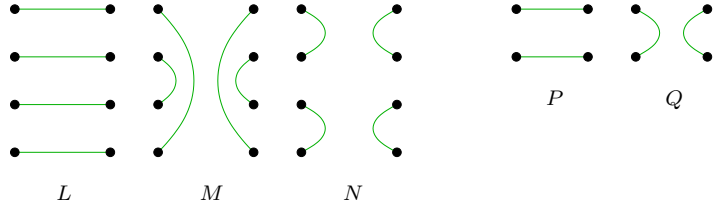}
\end{center}
\caption{The three even symmetric diagrams in of $\sTL^e_2$ and the
two Temperley-Lieb diagrams in $\TL_2$.} \label{fig:pfaf10}
\end{figure}

Let $L$, $M$, $N$ be the elements of $\sTL^e_2$ and $P$, $Q$ be the
elements of $\TL_2$ as shown in Figure \ref{fig:pfaf10}. Then
$\Pfaf_L(A) = yz-xt$, $\Pfaf_M(A) = xt-yz$, $\Pfaf_N(A) = xt$,
$\Imm_P(B) = xt-yz$, and $\Imm_Q(B) = yz$. One can check that the
coefficients of immanants in pfaffinants agree with Theorem
\ref{thm:dec}.
\end{example}

\subsection{Quadratic relations between TL-pfaffinants and
TL-immanants} Let $A$ be a skew-symmetric $2n \times 2n$ matrix. The
following formula is well known, see for example~\cite{St}.

\begin{theorem} \label{thm:square}
We have $\pf(A)^2 = {\rm det}(A)$.
\end{theorem}

More generally, we have the equation $\pf_{I, \bar I}(A)^2 =
\Delta_{I,I}(A) \Delta_{\bar I,\bar I}(A)$.  Let $S = I \cup (\hat
I)^\wedge$.  Applying Theorems~\ref{thm:decomp}
and~\ref{thm:immdecomp} we obtain the quadratic relationship
\begin{equation} \label{eq:pfafimm}
\left(\sum_{D \in \D(I)} \Pfaf_D(A) \right)^2 = \sum_{d \in \D(S)}
\Imm_d(B).
\end{equation}
It seems interesting to ask whether it is possible to refine
(\ref{eq:pfafimm}) to obtain simple quadratic relationships between
TL-pfaffinants and TL-immanants.  The simplest case arises from
Theorem~\ref{thm:square}: we have $\Pfaf_D(A)^2 = \Imm_d(A)$ where
$D \in \sTL^e_n$ and $d \in \TL_n$ are the respective diagrams
containing only horizontal edges (see Example~\ref{ex:pfafhor}). The
following table describes how to express TL-immanants of a
skew-symmetric $4 \times 4$ matrix in terms of the TL-pfaffinants.
Here we preserve the labeling of TL-pfaffinants as on Figure
\ref{fig:pfaf10}.

\begin{center}
\begin{tabular}{|c|c|}
\hline TL-Immanant & Sums of products of TL-pfaffinants
\\
\hline ${(4, 5), (3, 6), (2, 7), (1, 8)}$ & $\Pfaf_L^2$ \\
\hline ${(1, 2), (4, 5), (3, 6), (7, 8)}$ & $-\Pfaf_L^2$ \\
\hline ${(1, 2), (4, 5), (6, 7), (3, 8)}$ & $-\Pfaf_L \Pfaf_M$ \\
\hline ${(1, 2), (5, 6), (4, 7), (3, 8)}$ & $-\Pfaf_L^2 - \Pfaf_L \Pfaf_N$ \\
\hline ${(2, 3), (4, 5), (6, 7), (1, 8)}$ & $2 \Pfaf_L \Pfaf_M$ \\
\hline ${(2, 3), (1, 4), (6, 7), (5, 8)}$ & $\Pfaf_M^2$ \\
\hline ${(1, 2), (3, 4), (6, 7), (5, 8)}$ & $\Pfaf_L^2 + \Pfaf_L \Pfaf_M + \Pfaf_L \Pfaf_N + \Pfaf_M \Pfaf_N$ \\
\hline ${(1, 2), (3, 4), (5, 6), (7, 8)}$ & $2 \Pfaf_L^2 + 2 \Pfaf_L \Pfaf_N + \Pfaf_N^2$ \\
\hline
\end{tabular}
\end{center}

However for $n > 2$ the TL-immanants cannot be expressed in a
similar manner through TL-pfaffinants.  For example, with $n=3$ the
immanant corresponding to the diagram with edge set $\{(2, 3), (4,
5), (6, 7), (8, 9), (10, 11), (1, 12)\}$ does not lie in the span of
the products of the TL-pfaffinants. It remains unclear if any
relation between TL-immanants and TL-pfaffinants of a skew symmetric
matrix can be established in general.

\section{Schur $Q$-positivity} \label{pos}
In this section we discuss some conjectural applications of
TL-pfaffinants to positivity properties of Schur $Q$-functions. Many
of our results and conjectures can be stated alternatively in terms
of Schur $P$-functions, but we will not do so explicitly.

\subsection{Shifted tableaux }
For further details concerning the material of this section we refer
the reader to~\cite{Mac}.

Let $\lambda = \lambda_1 > \lambda_2 > \cdots > \lambda_l > 0$ be a
strict partition of integers.  We will not distinguish between
$\lambda$ and its {\it shifted diagram} $S(\lambda)$ obtained by
shifting the $i$-th row of the usual (Young) diagram $(i-1)$ squares
to the right, for each $i$.  More generally, if $\lambda$ and $\mu$
are two strict partitions so that $S(\mu) \subset S(\lambda)$ then
the skew shifted diagram is denoted $\lambda/\mu$.  Our notation for
diagrams follows the English notation, so that Young diagrams are
top-left justified.

A {\it shifted tableaux} $T$ with shape $\sh(T) = \lambda/\mu$ is a
filling of the shifted diagram $\lambda/\mu$ with the numbers
$1',1,2',2',\ldots$ so that
\begin{enumerate}
\item the rows and columns are weakly increasing under the order $1' <
1 < 2' < 2 < \ldots$
\item there is at most one occurrence of $i'$ in a row
\item there is at most one occurrence of $i$ in a column.
\end{enumerate}
The {\it weight} $\wt(T)$ of a shifted tableau is the composition
$\alpha = (\alpha_1,\alpha_2,\ldots)$ where $\alpha_i$ is equal to
the combined number of the letters $i$ and $i'$ used in $T$.  The
Schur $Q$-function $Q_{\lambda/\mu}(x)$ is defined as
$$
Q_{\lambda/\mu}(x) = \sum_{T\, : \, sh(T) = \lambda/\mu} x^{\wt(T)}
$$
where $x^\alpha = x_1^{\alpha_1} x_2^{\alpha_2} \cdots$.  Though it
is not immediate from the definition, the function
$Q_{\lambda/\mu}(x)$ is a symmetric function in the variables
$x_1,x_2,\ldots$.

\subsection{Schur $Q$-functions and pfaffians}
Schur $Q$-functions can be expressed as pfaffians, as follows. First
extend the notation of Schur $Q$-functions by defining $Q_{-r} = 0$
for $r > 0$ and $Q_{(r,s)} = -Q_{(s,r)}$.  Define the $l \times l$
skew symmetric matrix $A_{\lambda} = [Q_{(\lambda_i, \lambda_j)}]_{1
\leq i , j \leq l}$, where $\lambda_i$ is the $i$-th part of
$\lambda$, $l = l(\lambda)$ is the number of parts of $\lambda$.  By
possibly adding an extra zero part to $\lambda$, we may assume that
$l$ is even.  The following theorem can be found in \cite{Mac}.

\begin{theorem}
Let $\lambda$ be a strict partition.  Then $Q_{\lambda} =
\pf(A_{\lambda})$.
\end{theorem}
A skew version of this formula was proved by J\'{o}zefiak and
Pragacz~\cite{JP}.  Let $\lambda/\mu$ be a skew shifted shape where
$\lambda = \lambda_1 > \cdots \lambda_l > 0$ and $\mu = \mu_1 >
\mu_2
> \cdots > \mu_r \geq 0$.  We assume that $l + r$ is even.
Let $H = (h_{ij})$ be the $l \times r$ matrix with $h_{ij} =
Q_{\lambda_i - \mu_{r+1-j}}$. Define a skew symmetric matrix
$$A_{\lambda/\mu}= \left(
\begin{matrix}
A_\lambda & H\\
-H^t & 0
\end{matrix}
\right)
$$

%

We call the matrix $A_{\lambda/\mu}$ a {\it {$Q$-Jacobi-Trudi}}
matrix. If we allow in the definition $\lambda$ and $\mu$ to
possibly be non-strict partitions then we call $A_{\lambda/\mu}$ a
generalized $Q$-Jacobi-Trudi matrix.

\begin{theorem}[\cite{JP, St}] \label{thm:JP} Let $\lambda/\mu$ be a shifted skew
shape.  Then $Q_{\lambda/\mu} = \pf(A_{\lambda/\mu})$.
\end{theorem}
\begin{remark}
In~\cite{JP} the matrix Let $\tilde H = (\tilde h_{ij})$ be the $l
\times r$ matrix with $\tilde h_{ij} = Q_{\lambda_i - \mu_{j}}$ is
used.  Using $\tilde H$ to define $\tilde A_{\lambda/\mu}$, one then
has $\pf(A_{\lambda/\mu}) = (-1)^{r \choose 2} \pf(\tilde
A_{\lambda/\mu})$.

\end{remark}

\subsection{Schur $Q$-positivity and pfaffinants}
As we saw in Section~\ref{sec:npos}, network positivity of an
element $f \in \CP_n$ depends on the decomposition of $f$ into
diagram pfaffinants $\Pfaf'_D$.  Somewhat more surprisingly, we
conjecture that this decomposition is also related to Schur
$Q$-positivity.

\begin{conjecture} \label{conj:con1}
Suppose $f \in \CP_n$ can be expressed as $f = \sum c_D \Pfaf'_{D}$
with non-negative coefficients $c_D$.  Then for any generalized
$Q$-Jacobi-Trudi matrix $A_{\lambda/\mu}$, the evaluation
$f(A_{\lambda/\mu})$ is a nonnegative linear combination of Schur
$Q$-functions.
\end{conjecture}

Conjecture~\ref{conj:con1} parallels known Schur-positivity
properties of the TL-immanants $\Imm_d(A)$.  It is known~\cite{RS2}
that $\Imm_d(A)$ is nonnegative whenever $A$ arises from a planar
network and that the evaluations $\Imm_d(H_{\lambda,\mu})$ on
Jacobi-Trudi matrices $H_{\lambda,\mu}$ are nonnegative.

We can prove a weaker version of Conjecture~\ref{conj:con1}.

\begin{theorem} \label{thm:moncon1}
Suppose $f \in \CP_n$ can be expressed as $f = \sum c_D \Pfaf'_{D}$
with non-negative coefficients $c_D$.  Then for any shifted shape
$\lambda/\mu$, the evaluation $f(A_{\lambda/\mu})$ is a nonnegative
linear combination of monomial symmetric functions.
\end{theorem}

\begin{proof}
Stembridge has constructed a network $N_{\lambda/\mu}$ such that
$A_{\lambda/\mu} = A(N_{\lambda/\mu})$ (see \cite[Theorem 6.2]{St}).
By Theoerem~\ref{thm:JP}, $\pf(A_{\lambda/\mu}) = Q_{\lambda/\mu}$
also calculates the Schur $Q$-function.  It remains to note that for
each $D \in \sTL_n$ the evaluation $\hat
\Pfaf'_{D}(N_{\lambda/\mu})$ is a monomial positive formal power
series.  Finally, since $f \in \CP_n$ we know $f(A_{\lambda/\mu})$
must be a linear combination of Schur $Q$-functions and thus
symmetric.  We conclude that $f( A_{\lambda/\mu})$ is a nonnegative
linear combination of monomial symmetric functions.
\end{proof}

\subsection{Schur $Q$-functions and cell transfer}
In~\cite{LP}, we introduced the notion of a {\it $\T$-labeled
poset}. Let $Z$ be a totally ordered set (in~\cite{LP} we chose $Z =
\N$, but the results there generalize easily). Let $P$ be a poset
and $O$ be an assignment of a weakly increasing function $O(s
\lessdot t): Z \to Z \cup {\infty}$ to each cover relation $s
\lessdot t$ of $P$. A $(P,O)$-tableau $T$ is a function $T: P \to Z$
so that for each cover relation $s \lessdot t$ in $P$ we have $T(t)
\leq O(s,t)(T(s))$.

The boxes in a shifted diagram form a poset also denoted
$\lambda/\mu$.  The cover relations $s \lessdot t$ correspond to
boxes $s,t \in \lambda/\mu$ such that $s$ is immediately above or
immediately to the left of $t$.  Let $Z = \{1' < 1 < 2' < 2 <
\ldots\}$.  For a letter $z \in Z$ we denote by $z'$ the letter
obtained by either removing or adding a prime and for $z \in Z$ we
let $z \mapsto z+1$ denote the obvious operation which preserves
primes. Define the functions $f^r, f^c$ by
\begin{align*}
f^r(z) &= \begin{cases} z & \mbox{if $z$ is not primed,} \\
                        z' & \mbox{if $z$ is primed.} \end{cases} &
f^c(z) &= \begin{cases} z & \mbox{if $z$ is primed,} \\
                        z' + 1 & \mbox{if $z$ is primed.} \end{cases}
\end{align*}
Let $O_{\lambda/\mu}$ be the (edge) labeling of $\lambda/\mu$ such
that every cover relation along a row is labeled with $f^r$ and
every cover relation along a column is labeled with $f^c$.

The following result is immediate.
\begin{proposition}
Shifted tableaux of shape $\lambda/\mu$ are
$(\lambda/\mu,O_{\lambda/\mu})$-tableaux.
\end{proposition}

Now we define the operations $\wedge$ and $\vee$ on pairs of strict
partitions (see~\cite{LP,LPP}).  Namely, for partitions
$\lambda=(\lambda_1, \lambda_2, \dots,\lambda_n)$ and $\mu=(\mu_1,
\mu_2, \dots,\mu_n)$, define $\lambda\vee
\mu:=(\max(\lambda_1,\mu_1),\max(\lambda_2,\mu_2),\dots)$ and
$\lambda\wedge
\mu:=(\min(\lambda_1,\mu_1),\min(\lambda_2,\mu_2),\dots)$.  We may
have to add trailing zeroes before applying these operations.  These
operations send pairs of strict partitions to strict partitions. The
definition can be extended to skew shifted diagrams as follows:
$(\lambda/\mu)\vee(\nu/\rho) := (\lambda\vee\nu)/(\mu\vee\rho)$ and
$(\lambda/\mu)\wedge(\nu/\rho) :=
(\lambda\wedge\nu)/(\mu\wedge\rho)$.

In~\cite{LP} it was shown that the operations $\wedge$ and $\vee$
(called {\it cell transfer operations}) for $\T$-labeled posets give
rise to a range of monomial positivity results. By a slight
modification of~\cite[Theorem 3.6]{LP} one obtains a bijective proof
that $Q_{(\lambda/\mu)\vee(\nu/\rho)}
Q_{(\lambda/\mu)\wedge(\nu/\rho)} - Q_{\lambda/\mu}Q_{\nu/\rho}$ is
monomial positive.

The proof of the following stronger result will appear elsewhere.

\begin{theorem}\label{thm:peak}
Let $\lambda/\mu$ and $\nu/\rho$ be skew shifted shapes.  Then the
difference $$Q_{(\lambda/\mu)\vee(\nu/\rho)}
Q_{(\lambda/\mu)\wedge(\nu/\rho)} - Q_{\lambda/\mu}Q_{\nu/\rho}$$ is
a nonnegative sum of Stembridge's peak functions $K_\alpha$.
\end{theorem}

We will not give the definition of the peak functions $K_\alpha$
here and refer the reader to~\cite{Ste2} for full details.  The
$K_\alpha$ form a basis for a subalgebra $\P$ of the algebra
quasi-symmetric functions and the $K_\alpha$ take the place of the
fundamental quasi-symmetric functions in $\P$.  The Schur
$Q$-functions $Q_{\lambda/\mu}$ lie in this subalgebra $\P$ and are
known to be positive in the basis $\{K_\alpha\}$.  We now make the
following stronger conjecture.

%

\begin{conjecture} \label{conj:con2}
Let $\lambda/\mu$ and $\nu/\rho$ be skew shifted shapes.  Then the
difference $$Q_{(\lambda/\mu)\vee(\nu/\rho)}
Q_{(\lambda/\mu)\wedge(\nu/\rho)} - Q_{\lambda/\mu}Q_{\nu/\rho}$$ is
a non-negative combination of Schur $Q$-functions.
\end{conjecture}

Conjecture~\ref{conj:con2} is a Schur $Q$-function version of what
we call the {\it cell transfer} theorem.  The monomial positivity
version was proved in~\cite{LP}, the fundamental quasi-symmetric
function version in~\cite{LP2} and the Schur positivity version
in~\cite{LPP}.  As explained in the introduction of~\cite{LP2},
these positivity phenomena arise from a collection of data: (a) a
class of posets, (b) a ring containing the generating functions of
``tableaux'', (c) a basis of this ring, and (d) a set of skew
functions. In our case, (a) the posets are shifted Young diagrams,
(b) the ring is the subalgebra of the ring of symmetric functions
generated by the odd power sums, (c) the basis is the set of Schur
$Q$-functions for non-skew shifted shapes, and (d) the skew
functions are the Schur $Q$-functions labeled by skew shifted
shapes.

\begin{proposition} \label{prop:con2}
Conjecture \ref{conj:con1} implies Conjecture \ref{conj:con2} when $\mu = \rho =
\emptyset$.
\end{proposition}

\begin{proof}
Let $\pi$ be the (possibly no longer strict) partition obtained from
taking the union of the parts of $\lambda$ and $\nu$.  While $\pi$
is not necessarily a strict partition, we can still formally define
the matrix $A_{\pi}$ as above. Clearly, $\pf_{I, \bar I}(A_{\pi}) =
Q_{\lambda}Q_{\nu}$ for the appropriate choice of $I$. Now recall
the definition of $\min(I,\bar I)$ from Section~\ref{sec:comppfaf}.
We have $$\pf_{\min(I, \bar I),\overline{\min(I,\bar I)}}(A_{\pi}) -
\pf_{I,\bar I}(A_{\pi}) = Q_{(\lambda/\mu)\vee(\nu/\rho)}
Q_{(\lambda/\mu)\wedge(\nu/\rho)} - Q_{\lambda/\mu}Q_{\nu/\rho}.$$
By the proof of Proposition~\ref{prop:comppfaf}, the difference
$\pf_{\min(I, \bar I),\overline{\min(I,\bar I)}} - \pf_{I,\bar I}$
is a nonnegative linear combination of the TL-pfaffinants $\Pfaf_D$.
Conjecture~\ref{conj:con1} implies that $\Pfaf_D(A_{\pi})$ is Schur
$Q$-positive, from which the result follows.
\end{proof}

\subsection{Further Schur $Q$-positivity conjectures}
The usual Schur function analogue of Conjecture~\ref{conj:con2} was
established in~\cite{LPP}.
\begin{theorem}[\cite{LPP}] \label{thm:LPP}
Let $\lambda/\mu$ and $\nu/\rho$ be shifted shapes.  Then the
difference $s_{(\lambda/\mu)\vee(\nu/\rho)}
s_{(\lambda/\mu)\wedge(\nu/\rho)} - s_{\lambda/\mu}s_{\nu/\rho}$ is
a non-negative combination of Schur functions.
\end{theorem}

Theorem~\ref{thm:LPP} was used to resolve a number of conjectures of
Fomin, Fulton, Li, Poon~\cite{FFLP}, of Lascoux, Leclerc,
Thibon~\cite{LLT} and of Okounkov~\cite{Oko}.  We now state the
shifted analogue of the Fomin-Fulton-Li-Poon conjecture.

For two partitions $\lambda$ and $\mu$, let $\lambda \cup
\mu=(\nu_1,\nu_2,\nu_3,\dots)$ be the partition obtained by
rearranging all parts of $\lambda$ and $\mu$ in the weakly
decreasing order. Let $\mathrm{sort}_1(\lambda,\mu) :=
(\nu_1,\nu_3,\nu_5,\dots)$ and $\mathrm{sort}_2(\lambda,\mu) :=
(\nu_2,\nu_4,\nu_6,\dots)$. It is not hard to see that if $\lambda$
and $\mu$ are strict then so are $\mathrm{sort}_1(\lambda,\mu)$ and
$\mathrm{sort}_2(\lambda,\mu)$.

\begin{conjecture} \label{con3}
Let $\lambda, \mu$ be two shifted shapes.  Then
$Q_{\mathrm{sort}_1(\lambda,\mu)} Q_{\mathrm{sort}_2(\lambda,\mu)} -
Q_{\lambda} Q_{\mu}$ is a nonnegative linear combination of Schur
$Q$-functions.
\end{conjecture}

\begin{prop}
\label{prop:QFFLP} Conjecture \ref{conj:con2} implies Conjecture
\ref{con3}.
\end{prop}

\begin{proof}
First note that if $\lambda/\mu$ and $\nu/\rho$ are skew shifted
diagrams obtained from each other via a translation then
$Q_{\lambda/\mu} = Q_{\nu/\rho}$. For a shifted shape $\lambda$, let
$\lambda_\downarrow$ denote the skew shifted shape obtained by
translating $\lambda$ down one row (and hence also one step to the
right). We will assume that $\lambda_\downarrow$ is presented as
$\nu/\rho$ where $\nu_1 = \rho_1$ is very large (much larger than
any other parts involved in the proof).  If $\nu/\rho$ is a shifted
shape so that $\nu_1 = \rho_1$ we let $(\nu/\rho)_{\uparrow}$ denote
the shifted shape obtained by translating one row up (and hence also
one step to the left).

We can construct $\pi = \mathrm{sort}_1(\lambda,\mu)$ and $\theta =
\mathrm{sort}_2(\lambda,\mu)$ from $\lambda$ and $\mu$ by a sequence
of the operations $\vee$ and $\wedge$.  Suppose that we have
(strict) partitions $(\rho,\nu)$ so that $\rho \cup \nu = \lambda
\cup \mu$. Let us suppose that $\rho$ agrees with $\pi$ up to the
$i$-th part and that $\nu$ agrees with $\theta$ up to the $i$-th
part. If $\rho_{i+1} \neq \pi_{i+1}$, then it must be the case that
$\nu_{i+1} = \pi_{i+1}$.  In this case we replace $(\rho,\nu)$ by
$(\rho^*,\nu^*) = (\rho \vee \nu, \rho \wedge \nu)$. One checks that
$\rho^*$ agrees with $\pi$ up to the $(i+1)$-th part and $\nu^*$
agrees with $\theta$ up to the $i$-th part.  If $(\nu^*)_{i+1} \neq
\theta_{i+1}$ then $(\rho^*)_{i+2} = \theta_{i+1}$.  We now replace
$(\rho^*,\nu^*)$ by $(\rho^{**},\nu^{**}) = (\rho^* \wedge
(\nu^*)_\downarrow, (\rho^* \vee (\nu^*)_\downarrow)_\uparrow)$. One
checks that $\rho^{**}$ still agrees with $\pi$ up to the $(i+1)$-th
part and $\nu^{**}$ now agrees with $\theta$ up to the $(i+1)$-th
part.  After a finite number of iterations of the map $(\rho,\nu)
\mapsto (\rho^{**},\nu^{**})$ applied to $(\lambda,\mu)$, one
obtains $(\pi,\theta)$.

If we apply Conjecture~\ref{conj:con2} to the Schur $Q$-functions
indexed by the pairs of partitions $(\rho,\nu)$ we see that for each
iteration of the above map $ Q_{\rho^{**}}Q_{\nu^{**}}-Q_\rho Q_\nu$
is Schur $Q$-positive.  This proves the theorem.
\end{proof}

Our proof here is very similar to an analogous proof in~\cite{LPP},
where left and right shifts are used instead of our up and down
translations.  It would be interesting to generalize other Schur
positivity results and conjectures to the shifted case.

We note the following result, which follows from
Theorem~\ref{thm:peak} and the proof of
Proposition~\ref{prop:QFFLP}.

\begin{prop}
Let $\lambda, \mu$ be two shifted shapes.  Then
$Q_{\mathrm{sort}_1(\lambda,\mu)} Q_{\mathrm{sort}_2(\lambda,\mu)} -
Q_{\lambda} Q_{\mu}$ is a nonnegative linear combination of peak
functions.
\end{prop}


\section{Proof of Theorem~\ref{thm:welldefined}}
\label{sec:Rei}  Let $A$ and $B$ denote two nice embeddings of
$\nu(\pi)$ and denote by $f_D(A)$ and $f_D(B)$ the weight generating
function of uncrossings defined by $A$ and $B$ respectively (see
Section~\ref{sect:dpfaf}).  By replacing $A$ or $B$ with a small
deformation which is combinatorially equivalent we may assume even
if we draw all the edges of $A$ and $B$ that (a) no two edges have a
point of tangency and (b) no three strings cross at a single point.
However, an edge of $A$ and an edge of $B$ may intersect more than
once.

\begin{figure}[ht]
\begin{center}
\input{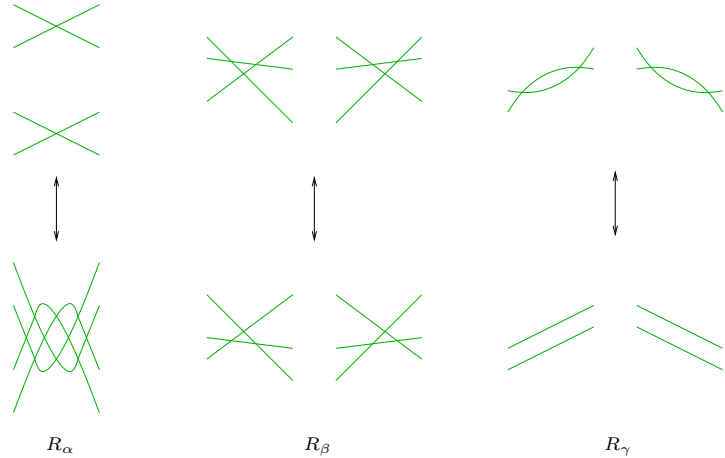}
\end{center}
\caption{The three types of Reidemeister-like moves
$R_\alpha,R_\beta, R_\gamma$.} \label{pfaf1}
\end{figure}

We now argue that $A$ and $B$ are connected by a sequence of three
types of Reidemeister-like moves, denoted $R_\alpha$, $R_\beta$ and
$R_\gamma$, as shown in Figure \ref{pfaf1}.  Let $(i,j)$ be an edge
in $\nu(\pi)$.  To change $A$ to $B$, we move the embedding of
$(i,j)$ in $A$ continuously until it agrees with the embedding of
$(i,j)$ in $B$; and we repeat for each edge of $\nu(\pi)$.  Note
that we will always move the mirror symmetric edge simultaneously so
that the diagram is always mirror symmetric.  There are three types
of ``singularities'' which may occur during this process, changing
the combinatorial type of the embedding.  These singularities
violate the conditions (a) and (b) above.

\begin{itemize}
\item[$R_\alpha$:]
If the singularity occurs on the vertical axis of symmetry then one
obtains a quadruple intersection between two pairs of mirror
symmetric edges, violating both conditions (a) and (b).  The
Reidemeister move $R_\alpha$ allows one to pass from one side of the
singularity to the other.

\item[$R_\beta$:]
If the singularity is a paired singularity, it may involve three
edges crossing at the same point, giving the move $R_\beta$.

\item[$R_\gamma$:]
If the singularity is a paired singularity, it may involve a point
of tangency, giving the move $R_\gamma$.
\end{itemize}

Note that $R_\alpha$ allows us to permute the crossing points on the
vertical axis of symmetry, while $R_\beta$ and $R_\gamma$ allow us
to do all the other required changes.  During this process the rule
that no two edges crossing more than once can be violated (by moves
$R_\alpha$ or $R_\gamma$).

To complete the proof we show that $f_D(A) = f_D(A')$ if $A$ and
$A'$ are related by a Reidemeister-like move.

\begin{itemize}
\item[$R_\alpha$]: We may use the move $R_\gamma$ (preserving $f_D$) to replace
the initial and final pictures with the two intermediate ones shown
in Figure~\ref{pfaf14}. Now using the the calculation of
Example~\ref{ex:1423} we may obtain the one intermediate picture
from the other while again preserving $f_D$.

\begin{figure}[ht]
\begin{center}
\input{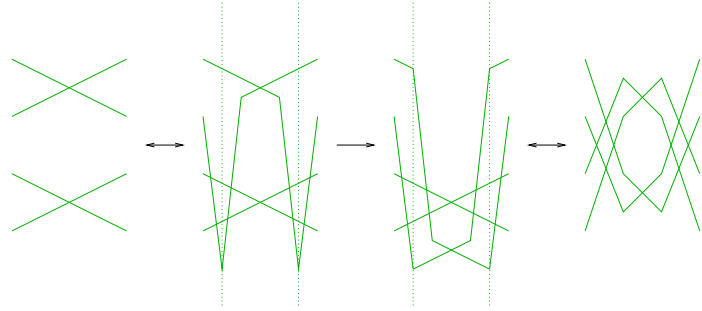}
\end{center}
\caption{Verification of the $R_{\alpha}$ move.} \label{pfaf14}
\end{figure}

\item[$R_\beta$]:  There are three pairs of (mirror-symmetric)
crossings, giving a total of 8 uncrossings for the initial and final
pictures.  Denote the three edges coming from the left by $a,b,c$
from top to bottom and the three edges exiting to the right by
$a',b',c'$.  An uncrossing of this local picture will give a
matching of $a,a',b,b',c,c'$ together with a weight.  One obtains
the following table for the weights of the 8 uncrossings, showing
that the weight generating functions agree for each matching. The
``initial'' embedding here is the top picture in Figure~\ref{pfaf1}.

\begin{center}
\begin{tabular}{|c|c|c|}
\hline Matching & Initial embedding & Final embedding
\\
\hline $\{(a,a'),(b,b'),(c,c')\}$ & 1 & 1\\
\hline $\{(a,b),(a',b'),(c,c')\}$ & $1 + 1 + 1 - 2 = 1$&  1\\
\hline $\{(a,a'),(b,c),(b',c')\}$ & 1 & $1 + 1 + 1 - 2 = 1$\\
\hline $\{(a,c'),(a',b'),(b,c)\}$ & -1& -1\\
\hline $\{(a,b),(c,a'),(b',c')\}$ & -1 & -1 \\
\hline
\end{tabular}
\end{center}

\item[$R_\gamma$:]  For the initial (top) embedding, the picture has
4 uncrossings.  Three of these 4 uncrossings give a matching (the
vertical one) which does not occur for the final (bottom) embedding,
but their weights (respectively $2$,$-1$,$-1$) cancel out.  For the
other (horizontal) matching we obtain the same contribution of 1 for
both the initial and final embeddings.

\end{itemize}

This completes the proof of Theorem~\ref{thm:welldefined}.

\end{document}